\documentclass[11pt, reqno]{amsart}
\usepackage{amsmath, amsfonts, amssymb, amscd, amsthm,mathrsfs, enumerate, enumitem, color, xcolor, cancel, mathtools}
\usepackage{ hyperref, times}
\usepackage[all]{xy}

\theoremstyle{plain}
\newtheorem{theo}{Theorem}[section]
\newtheorem{lem}[theo]{Lemma}
\newtheorem{prop}[theo]{Proposition}

\newtheorem{cor}[theo]{Corollary}
\theoremstyle{definition}
\newtheorem{rem}[theo]{Remark}

\newtheorem{definition}[theo]{Definition}
\newenvironment{pf}{\noindent{\it Proof. }}{$\square$\par\medskip}
\newenvironment{pfns}{\noindent{\it Proof. }}{\par\medskip}

\theoremstyle{plain}

\newtheorem{theorem}[theo]{Theorem}

\theoremstyle{definition}

\renewcommand{\=}{:=}

\newcommand{\rank}{\operatorname{rank}}
\newcommand{\beq}{\begin{equation}}
\newcommand{\eeq}{\end{equation}}
\renewcommand{\a}{\alpha}

\renewcommand{\d}{\delta}

\newcommand{\f}{\varphi}

\renewcommand{\l}{\lambda}

\renewcommand{\r}{\rho}

\renewcommand{\t}{\tau}

\renewcommand{\L}{\Lambda}

\newcommand{\bC}{\mathbb{C}}

\newcommand{\bK}{\mathbb{K}}

\newcommand{\bR}{\mathbb{R}}
\newcommand{\bZ}{\mathbb{Z}}

\newcommand{\bN}{\mathbb{N}}

\newcommand{\bX}{\mathbb{X}}

\renewcommand{\gg}{\mathfrak{g}}

\newcommand{\gi}{\mathfrak{i}}

\newcommand{\gm}{\mathfrak{m}}
\newcommand{\gn}{\mathfrak{n}}

\newcommand{\gp}{\mathfrak{p}}

\newcommand{\gD}{\mathfrak{D}}

\newcommand{\gX}{\mathfrak{X}}

\newcommand{\gU}{\mathfrak{U}}

\newcommand{\ggl}{\mathfrak{gl}}

\newcommand\SU{\mathrm{SU}}
\newcommand\U{\mathrm{U}}

\newcommand{\cB}{\mathcal{B}}
\newcommand{\cC}{\mathcal{C}}
\newcommand{\cD}{\mathcal{D}}

\newcommand{\cH}{\mathcal{H}}

\newcommand{\cK}{\mathcal{K}}

\newcommand{\Jst}{J_{\operatorname{st}}}

\renewcommand{\square}{\kern1pt\vbox
{\hrule height 0.6pt\hbox{\vrule width 0.6pt\hskip 3pt
\vbox{\vskip 6pt}\hskip 3pt\vrule width 0.6pt}\hrule height0.6pt}\kern1pt}

\DeclareMathOperator\Aut{Aut\;}
\DeclareMathOperator\aut{aut}

\DeclareMathOperator\ad{ad}

\DeclareMathOperator{\Span}{Span}

\newcommand\Hom{\operatorname{Hom}}

\newcommand{\wt}{\widetilde}
\newcommand{\wh}{\widehat}

\newcommand{\be}{\begin{equation}}
\newcommand{\ee}{\end{equation}}

\def\<#1,#2>{\langle\,#1,\,#2\,\rangle}
\newcommand{\arr}{\begin{array}{rlll}}
\newcommand{\ea}{\end{array}}
\newcommand{\bea}{\begin{eqnarray}}
\newcommand{\eea}{\end{eqnarray}}
\newcommand{\bean}{\begin{eqnarray*}}
\newcommand{\eean}{\end{eqnarray*}}

\def\sideremark#1{\ifvmode\leavevmode\fi\vadjust{
\vbox to0pt{\hbox to 0pt{\hskip\hsize\hskip1em
\vbox{\hsize3cm\tiny\raggedright\pretolerance10000
\noindent #1\hfill}\hss}\vbox to8pt{\vfil}\vss}}}

\newcounter{ssig}
\newcounter{ttig}
\newcommand{\under}[1]{{\underline{#1}\,}}

\newcommand{\degree}{\operatorname{degree}}
\newcommand{\type}{\operatorname{type}}

\newcommand{\bdot}{{\boldsymbol{\cdot}}}

\subjclass[2010]{32V05, 32V40,  22F30, 22F50, 57S25}
\keywords{Totally nondegenerate CR manifold; Beloshapka conjecture;  Tanaka structure; Tanaka symbol}

\thanks{{\it Acknowledgments}.  The research of the first author was supported in part by the grant from IPM, no. 96510425.
The second author was partially supported by the Project MIUR ``Real and Complex Manifolds: Geometry, Topology and  Harmonic Analysis'' and by GNSAGA of INdAM}

\begin{document}

\begin{abstract}   A CR manifold  $M$, with CR distribution $\cD^{10} \subset T^\bC M$,  is called {\it totally nondegenerate of depth $\mu$} if:
(a) the complex tangent space $T^\bC M$ is generated by all complex vector fields
that might be determined    by  iterated  Lie brackets between at most   $\mu$  fields  in $\cD^{10} + \overline{\cD^{10}}$;
(b) for each integer $2 \leq k \leq \mu-1$, the families of all vector fields that might be determined   by  iterated  Lie brackets between at most   $k$  fields  in $\cD^{10} + \overline{\cD^{10}}$  generate  regular complex distributions;
 (c) the ranks of the distributions  in (b)  have the {\it  maximal values}  that can be obtained amongst all CR manifolds
of the same CR dimension and satisfying (a) and (b)  -- this  maximality property is   the  {\it total nondegeneracy} condition.
In this paper, we  prove that, for any Tanaka symbol $\gm =  \gm^{-\mu} + \ldots + \gm^{-1}$ of a totally nondegenerate CR manifold of depth $\mu \geq 4$,
the full Tanaka prolongation of $\gm$ has trivial subspaces  of degree $k \geq 1$, i.e. it  has  the form  $\gm^{-\mu} + \ldots + \gm^{-1} + \gg^0$.  This result has various consequences.  For instance  it implies that any  (local) CR automorphism  of  a regular totally nondegenerate  CR manifold  is uniquely determined by its first order jet at a fixed point of  the manifold.  It also gives a complete proof of  a conjecture  by  Beloshapka on the group of automorphisms  of homogeneous totally nondegenerate CR manifolds.\\[-1cm]
\end{abstract}

\title[On the geometric order of
 totally nondegenerate CR manifolds
]{
On the geometric order of \\
totally nondegenerate CR manifolds
}

\author{Masoud Sabzevari}
\address{Department of Mathematics,
Shahrekord University, 88186-34141 Shahrekord, IRAN and School of
Mathematics, Institute for Research in Fundamental Sciences (IPM), P.
O. Box: 19395-5746, Tehran, IRAN}
\email{sabzevari@ipm.ir}

\author{Andrea Spiro}
\address{Scuola di Scienze e Tecnologie, Universit\`{a} di Camerino,
Via Madonna delle Carceri, 62032 Camerino, Macerata, ITALY}
\email{andrea.spiro@unicam.it}


\maketitle
\section{Introduction}
An {\it (abstract) CR manifold} is a real manifold $M$, equipped with a $\cC^\infty$ complex distribution $\cD^{10} \subset T^\bC M$ such that: (i) $\cD^{10} \cap \overline{\cD^{10}} = \{0\}$;  (ii)  for all vector fields $X, Y \in \cD^{10}$ the Lie bracket $[X, Y]$ is also in $\cD^{10}$.   The rank of the distribution  $\cD^{10}$ is called {\it CR dimension}. The most  natural and  studied examples  are the   {\it embedded CR manifolds}, which are   the  smooth real submanifolds  $M \subset \bC^N$, $N \geq 2$,   satisfying  appropriate constant rank conditions that guarantee that  the family of all complex vector fields of  the holomorphic distribution of $\bC^N$  with  real and imaginary parts tangent to $M$,  generate a   complex distribution    $\cD^{10} \subset  T^\bC M$ of constant rank.  \par
\smallskip
In this paper, we study the geometric structures of   the  {\it totally nondegenerate CR manifolds of depth $\mu$}, a class of CR manifolds  introduced by Beloshapka in \cite{Be1}.  For motivation and main reasons of interest for such important class of CR structures we refer to the original paper  (see also \cite{Be2, Sa1, Sa2}). The properties which characterise  a totally nondegenerate CR manifold $(M, \cD^{10})$  of depth $\mu$ can be shortly  described as follows (see \S 3 for the detailed definition):
\begin{itemize}[itemsep=5pt, leftmargin=20pt]
\item[a)] the complexified tangent bundle $T^\bC M$ is spanned  by the vector fields in  $\cD^{10} + \overline{\cD^{10}}$  and  by all  iterated  Lie brackets determined by sets of  at most   $\mu$   vector  fields in $\cD^{10} + \overline{\cD^{10}}$;
\item[b)] for each  $2 \leq k \leq \mu-1$, the  vector fields in  $\cD^{10} + \overline{\cD^{10}}$ and   all  iterated  Lie brackets determined by sets of  at most    $k$
 vector fields  in $\cD^{10} + \overline{\cD^{10}}$ span a   regular complex distribution, denoted by $(\cD^k)^\bC$;
\item[c)] the ranks of the distributions $(\cD^k)^\bC$,  $2 \leq k \leq \mu-1$,   are  the  maximal possible ones  that can occur  in  any other  CR manifold of the same CR dimension and satisfying (a) and (b) -- this  maximality property is  called   {\it total  nondegeneracy} condition.
\end{itemize}
The  CR manifolds $M$   that  are appropriately osculated at each point $x \in M$ by  a fixed homogeneous CR manifold  satisfying (a) -- (c) are called {\it regular}  (or {\it of uniform type}).\par
\smallskip
To be rigorous,   we have to  mention that in \cite{Be1}  Beloshapka  did not really  give the above definition, but actually  introduced  the   notion   of   {\it  germ of a totally nondegenerate  embedded CR submanifold} of $\bC^N$ and   presented it in terms of some appropriate  normal forms for its  defining equations. Nonetheless,  one can  directly  check that, if  one translates  everything  into the language of  germs of  embedded submanifolds,    our definition becomes completely  equivalent to  Beloshapka's  one.\par
\smallskip
According to the results by  Tanaka on the  equivalence problem for non-integrable distributions and for the CR structures (\cite{Ta}; see also \cite{AS}),
any totally nondegenerate CR manifold $(M, \cD^{10})$ of depth $\mu$   is  canonically associated with a family of pairs $(\gm_x, J_x)$, one  per each $x \in M$, formed by:
\begin{itemize}[itemsep=5pt, leftmargin=20pt]
\item[--] a negatively graded Lie algebra $\gm_x = \gm^{-\mu}_x + \gm^{-\mu +1}_x + \ldots + \gm^{-1}_x$,  called {\it Tanaka's symbol  at $x$},  with  $\dim \gm_x = T_x M$ and  with $(\gm^{-1}_x)^\bC \simeq \cD^{10}_x \oplus \overline{\cD^{10}_x}$;
\item[--]  a complex  structure $J_x: \gm^{-1}_x \to \gm_x^{-1}$,  whose $(\pm i)$-eigenspaces $\gm^{10}_x$ and $\gm^{01}_x$ in the complexified space $(\gm^{-1})^\bC \simeq  \cD^{10}_x \oplus \overline{\cD^{10}_x}$ coincide with   the spaces  $\cD^{10}_x$ and  $\overline{\cD^{10}_x}$, respectively (see \S 2  for details).
\end{itemize}
A CR manifold is  regular if and only if all  pairs $(\gm_x, J_x)$ are  isomorphic to  the same one,  say  $(\gm, J)$.  Hence, if  the CR manifold  is regular,  we may say that at each point $x$ it is ``osculated'' by  the CR manifold, given  by the  simply connected Lie group $G^\gm$  with $Lie(G^\gm) = \gm$ equipped  with the  unique   $G^\gm$-invariant  distribution $\cD^{\gm 10} \subset T^\bC G^\gm$  with  $\cD^{\gm 10}|_e  = \gm^{10}$. \par
\smallskip
We remark that, in \cite{Be1},  Beloshapka considered a particular class of embedded totally nondegenerate  CR manifolds, called {\it model surfaces},    and proved  that all of them are acted on transitively  by a group of CR transformations  whose  Lie algebra  is isomorphic to  one of the above described graded Lie algebras  $Lie(G^\gm) = \gm$. This immediately implies   that, up to  coverings,  Beloshapka's model surfaces are nothing but
embedded realisations in   $\bC^N$ of the  above described osculating manifolds $G^\gm$.\par
\par
\smallskip
This paper is devoted to the proof of the following  (Theorem \ref{main}): {\it if  $\mu \geq 4$,   the Lie algebra of the group $\Aut(G^\gm)$ of  CR automorphisms of the homogenous  CR manifold $G^\gm$ of depth $\mu$    is the graded Lie algebra
$$Lie(\Aut(G^\gm)) = \gm + \gg^0 = \gm^{-\mu} + \ldots + \gm^{-1} + \gg^0\ ,\eqno{(\ast)}$$
where $\gg^0$ is  the Lie algebra of the group of  automorphisms of  the graded Lie algebra $\gm$ that  leave  invariant the complex structure $J: \gm^{-1} \to \gm^{-1}$.} Note that if $n$ is the CR dimension, $\gg^0$ is naturally identifiable with a subalgebra of $ \ggl_n(\bC)$.  \par
\smallskip
Our main result naturally includes and completes  all previous results on the automorphisms of totally nondegenerate CR manifolds of depth $\mu \geq 4$, obtained in diverse papers by  Kossovski\v\i\ and  by the first author  (\cite{Kos, Sa1, Sa2}; see also \cite{SHAM1, SHAM2}). Indeed, our proof has been built  on several  ideas of the first author  and it can be considered as a realisation of a project,  which has been  outlined in his previous   papers. \par
\smallskip
By   Tanaka's theory (\cite{Ta,AS}), our result  has immediate interesting  consequences on the geometry of {\it any} regular totally nondegenerate CR manifold. For instance, it implies that for  any $n$ CR dimensional manifold  $M$ of this kind, the  dimension of the automorphism group is  less than or equal to $\dim M +  2 n^2$ and   any infinitesimal CR transformation is  completely determined by its first order jet at some point. An expanded discussion of such  geometric outcomes will be given in a future paper.\par
\smallskip
We also remark that, by  the existence of CR equivalences (up to coverings) between the   homogeneous CR  manifolds  $G^\gm$ and Beloshapka's  model surfaces, our result  gives also a complete proof for the cases $\mu \geq  4$  of a conjectured property, called `maximum conjecture'' in \cite{Be2}. There, Beloshapka wrote that  the Lie algebra of the automorphism groups of {\it any}  model  surface of depth $\mu \geq 3$ was expected  to be  precisely of the form  ($\ast$).
Since   this property   has been proved  for    $\mu = 3$   by Gammel' and Kossovski\v \i\ in \cite{GK},  by the  results of this paper  we may  claim that  Beloshapka's conjecture  is now confirmed   in all  cases. \par
\medskip
The paper is structured as follows. In \S 2 and \S 3, we give a  short review  of all  elements of Tanaka's theory of fundamental Lie algebras, prolongations, etc. , which we need for our proof, and  we introduce the notion of {\it universal fundamental CR algebra}. This can be considered as a CR analogue of the notion of  free Lie algebra (also called {\it universal  fundamental Lie algebra} in \cite{Ta}) and  allows a very convenient   characterisation  of Tanaka's symbols  of  totally nondegenerate CR structures.   Section \S 3 ends with a detailed definition of regular totally nondegenerate CR manifolds and  the  statement of our main result.  In \S 4 we give a crucial  result on linear independent sets   in universal fundamental CR algebras and in totally nondegenerate fundamental Lie algebras.  The proof of the   main theorem is given in  \S 5.\\[4pt]
{\it Acknowledgement.}  After posting our paper on  arXiv,  we  noticed that,
independently and  practically simultaneously to us,
our main result has been proven also  by Jan Gregorovi\v c  in \cite{Gr}.
The authors are  grateful to Ilya  Kossovski\v\i\ for pointing this preprint to us.
We are also  very grateful to Jo\"el Merker and the referee
for careful readings and truly helpful and constructive  remarks.
\bigskip
\section{Preliminaries}
\subsection{Fundamental algebras, CR Tanaka algebras and  CR manifolds}
\subsubsection{Fundamental algebras and  CR Tanaka algebras}
We recall that a {\it fundamental algebra}  of depth $\mu$ is a  negatively graded Lie algebra over $\bK = \bR$ or $\bC$ of the form
$$
\gm:=\gm^{-\mu} + \gm^{-\mu + 1} + \ldots + \gm^{-1}$$
with
$[\gm^{-i}, \gm^{-1}]  = \gm^{-i -1}$  for each  $i$.  Further, denoting by $\aut(\gm) = Lie(\Aut(\gm))$ the   Lie algebra
of the group  $\Aut(\gm)$  of  all automorphisms of the graded Lie algebra $\gm$ and   given a  subalgebra $\gg^0 \subset \aut(\gm)$,  the {\it Tanaka algebra} with isotropy $\gg^0$ is  the  non-positively graded  Lie algebra
$$\gm + \gg^0 = \gm^{-\mu} + \ldots + \gm^{-1} + \gg^0\ .$$
\par
\smallskip
When $\gm$  is real and is  equipped with a complex structure  $J: \gm^{-1}\rightarrow\gm^{-1}$ on the subspace $\gm^{-1}$ that   satisfies the so-called {\it integrability condition}
\beq \label{integrability} [J X, JY] = [X, Y]\qquad \text{for all}\ X, Y \in \gm^{-1}\ , \eeq
the pair $(\gm, J)$ is called {\it fundamental CR algebra}. The Tanaka  algebra   $\gm + \gg^0$,  with $\gg^0$ equal to the Lie algebra  $\gg^0 = \aut(\gm, J) = Lie(\Aut(\gm, J))$ of the group $\Aut(\gm, J)$ of $J$-preserving automorphisms is called  {\it CR Tanaka algebra}  determined by   $(\gm, J)$. \par
\medskip
\subsubsection{Regular CR manifolds of given type and associated Tanaka algebras} \label{2.1.2} Let $M$ be a real $n$-dimensional manifold.   A (integrable) {\it CR structure} on $M$ is a pair $(\cD, J)$,  formed by a distribution $\cD \subset TM$ of  rank $2 k \leq n$ and a smooth family of complex structures $J_x: \cD_x \to \cD_x$, $x \in M$,  with the property  that the complex distribution $\cD^{10} \subset T^\bC M$  formed  by the $+i$-eigenspaces of the complex structures $J_x: \cD^\bC \to \cD^\bC$, is involutive. A  diffeomorphism $f: M \to M$ of a CR manifold is called {\it CR automorphism} if $f_*(\cD) \subset \cD$ and $f_*(J) = J$. A vector field $X \in \gX(M)$, whose (local) flow consists of $1$-parameter families of (local) CR automorphisms is called an  {\it infinitesimal (local) CR automorphism} of $(M, \cD, J)$. The space of all infinitesimal CR transformations  is known to be a Lie algebra and we denote it by $\aut(M, \cD, J)$. \par
\smallskip
We recall that the  most relevant examples of CR manifolds are   given by  real submanifolds  of  $\bC^N$. More precisely,   if we denote by $\Jst$ the standard complex structure of $\bC^n$,  any $\cC^\infty$ real submanifold $M \subset \bC^N$, for which   the family $\cD^o$ of tangent subspaces  $\cD^o_x  = \{\ v \in T_x M\ :\  \Jst v = \sqrt{-1} v \in T_x M\ \}$  is a regular distribution of constant rank $2k$, then  the pair  $(\cD^o, J\= \Jst|_{\cD})$ is   a  CR structure on $M$. It is called the   {\it induced CR structure} of the regularly embedded submanifold $M \subset \bC^N$.    \par
\smallskip
We also remind that  if $(\cD, J)$ is  a CR structure on $M$, the integer $k = \dim M - \rank \cD$ is called {\it CR codimension} and,  for  induced CR structures
of regularly embedded CR generic submanifolds, this is also equal to  the   real codimension of the submanifold.\par
\medskip
There exists a canonical relation between CR Tanaka algebras and CR manifolds satisfying appropriate regularity conditions.  Consider an $n$-dimensional CR manifold $(M, \cD, J)$. Further,  for each  (real or complex) distribution  $\cK$ on  $M$,  let us  adopt the notational convention of indicating by $\under \cK$ the class of all local smooth vector fields taking values  of  $\cK$.   We may now consider  the sequence of   spaces $\gD_{-j}$, $ 1 \leq j$ of vector fields, defined inductively by
\beq \label{recursive} \gD_{-1} \= \under \cD\ ,\qquad
\gD_{-(j+1)}:=\gD_{-j}+[\gD_{-1}, \gD_{-j}] \ .\eeq
If  there exists a sequence of constant rank distributions $\cD_{-j} \subset TM$  such that  $\gD_{-j} = \under \cD_{-j}$ for each $j \geq 2$,  we say that the distribution  $\cD$ is  {\it regular}  (\cite{Ta}). In this case,  for   each point $x \in M$, we may define
$$\gm^{-1}(x) \= \cD|_x\ ,\qquad \gm^{-(j + 1)}(x) \= \cD_{-(j+1)}|_x/\cD_{-j}|_x$$
and, denoting by  $\mu$  the first integer for which  $\gm^{-(\mu+ k)} = 0$ for all $k \geq 1$,  we  set
$$\gm(x) = \gm^{-\mu}(x)+ \cdots + \gm^{-1}(x)\ .$$
The usual Lie brackets between  vector fields induce natural  Lie brackets on the sum of quotient spaces on $\gm(x)$  and  makes it  a fundamental algebra of depth $\mu$. Further, as a consequence of the   integrability  of  $J$, one can check that the complex structure $J_x: \cD_x = \gm^{-1}(x) \to \cD_x = \gm^{-1}(x)$  makes    $(\gm(x), J_x)$ a fundamental CR algebra.\par
All this allows to consider the following
\begin{definition} Given a fundamental CR algebra $(\gm, J)$, a  CR manifold $(M, \cD, J)$  is called {\it regular  of type $(\gm, J)$}  if\\
  a) the distribution  $\cD$ is   regular,  \\
   b) all fundamental CR algebras $(\gm(x), J(x))$, $x \in M$, are  isomorphic to $\gm$,\\
    c) $\dim \gm = \dim M$.   \\
The  CR Tanaka algebra $\gm + \gg^0$,  given by $(\gm, J)$,  is   called   {\it   associated  to $(M, \cD, J)$}.
\end{definition}
\par
 Given a fundamental CR algebra $(\gm, J)$, let $G^{\gm}$ to be the unique (up to an isomorphism) simply connected Lie group
  with $\gm = Lie(G^\gm)$ and consider the  left invariant distribution $\cD^\gm$, with $\cD^\gm|_e = \gm^{-1} \subset T_e G^\gm$, and the
  left invariant family of complex structures $J^\gm|_g: \cD^\gm_g \to \cD^\gm_g$ for $g\in G^{\gm}$, with $J^\gm|_e = J$.  The pair $(\cD^\gm, J^\gm)$ is directly seen to be an integrable CR structure on $G^\gm$ (\cite[\S 10.4]{Ta}). The homogeneous CR manifold $(G^\gm, \cD^\gm, J^\gm)$ constructed in this way is called the {\it standard CR model associated with $(\gm, J)$}.\par
  \smallskip
  We stress the fact that, being  $(G^\gm, \cD^\gm, J^\gm)$  homogeneous, it has a natural structure of real-analytic manifold with a real-analytic   CR structure. So, by classical embedding results for analytic CR manifolds,  {\it any such standard CR model admits a   local  regular CR embedding  in  $\bC^N$,  with
  $N =\dim_\bR G^\gm - \frac{1}{2} \dim_\bR \cD^\gm_o$}.  {\it  The germ of such embedding is unique up to local biholomorphisms}   (\cite[Props. 10.1 - 10.2]{Ta}).\par
 \medskip
If $(M, \cD, J)$ is regular and of type $(\gm, J)$, then  the class of all $J$-equivariant graded   isomorphisms $\xi_x: \gm(x) \to \gm$, $x \in M$,
determines a  principal  bundle $\pi: P^0 \to M$ over $M$, with structure group $G^0 \= \Aut(\gm,J)$. This bundle  is canonical  in the sense that it is preserved by  any CR  automorphism.  More precisely, in case  $f: M \to M$ is a diffeomorphism in $\Aut(M, \cD, J)$, the corresponding push-forward map $f_*: TM \to TM$ induces, in a standard way,  a   Lie algebras isomorphism $f_*|_x: \gm(x)\to \gm(f(x))$ between the fundamental algebras at the  points $x$ and $f(x)$, for each $x \in M$.  This allows to consider  the map
$$ \wh f(\xi_x) \=  \xi_x \circ (f_*|_x)^{-1}: \gm(f(x)) \to \gm\ ,$$
which transforms each isomorphism  $\xi_x: \gm(x) \to \gm$ of the fiber $ P^0|_x$ into a corresponding isomorphism $\wh f(\xi_x): \gm(f(x)) \to \gm$ of the fiber $P^0|_{f(x)}$. This defines a principle bundle automorphism $\wh f: P^0 \to P^0$, which projects onto the diffeomorphism $f: M \to M$ and it is uniquely determined by such projection.  Such an automorphism  $\wh f$ is called  {\it the canonical  lift of $f \in \Aut(M, \cD, J)$ to $P^0$}. Taking in consideration flows, the canonical correspondence $f \to \wh f$ naturally induces a canonical correspondence between  the vector fields  $X \in \aut(M, \cD, J)$ and appropriate vector fields $\wh X$ on $P^0$, which we call {\it canonical lifts} of the  infinitesimal CR transformations.\par
\medskip
\subsection{Prolongations of CR Tanaka algebras and CR   geometries of finite order}
\subsubsection{Prolongations of  Tanaka algebras}
Given  a Tanaka algebra $\gm+\gg^0$,  its {\it (full) prolongation}  is the graded Lie algebra
$$
(\gm+\gg^0)^\infty:=\gm^{-\mu} + \ldots  + \gm^{-1}  + \gg^0+\gg^1+\gg^2+\ldots\ , $$
in which the positively graded spaces $\gg^\ell$ with $1 \leq \ell < \infty$ and the Lie  brackets are inductively  defined  as follows. Given an integer $\ell \geq 1$,  assume that all spaces $\gg^i$,  $0 \leq i \leq \ell-1$,  and all  brackets $[X, Y]$ between  homogeneous elements  $X \in  \sum_{i = 0}^{\ell-1} \gg^i$ and  homogeneous element $Y \in \gm$ have been defined. Then,    set
\begin{multline}
\label{defg1}
\gg^\ell \= \{\ X  \in {\textstyle \sum_{-p = -\mu}^{-1}} \Hom(\gm^{-p},  \wt  \gg^{-p + \ell})\ : \\
[X(Y), Z] + [Y, X(Z)] = X([Y, Z]) \quad
\text{for each}
\ \  Y , Z \in \gm \ \}\ ,
\end{multline}
where $\wt \gg^i$  stands for   $\gm^i$ in case $i \leq -1$,  and  for $\gg^i$ in case $i \geq 0$.
Define also brackets between  homogeneous elements of the form $X^\ell \in  \gg^\ell$, $Y^{-k} \in \gm^{-k} $ by
\beq \label{Lie2-bis} [X^\ell, Y^{-k}] \= X^\ell(Y^{-k}) \in \wt \gg^{-k + \ell}\ .\eeq
In \cite{Ta}, Tanaka proved that there exists a unique way to define a graded Lie algebra structure on the direct sum  $ (\gm+\gg^0)^\infty = \gm  + \gg^0+ \sum_{i = 1}^\infty \gg^i$,  whose Lie brackets coincide with  the original  brackets between pairs in  $\gm + \gg^0$ and is equal to  \eqref{Lie2-bis} for each $k$ and $\ell$.
Note that, by construction, {\it if there is $k_o$ such that $\gg^{k_o + 1} = \{0\}$, then $\gg^{k_o + \ell} = \{0\}$ for all $\ell \geq 1$}.
\medskip
\subsubsection{Tanaka's towers of CR manifolds and lifts of  CR automorphisms } \label{2.2.2}
Let $(M, \cD, J)$ be a regular CR manifold of type $(\gm, J)$. According  to Tanaka's theory of differential systems (\cite{Ta, AS}),  there  exists a   sequence  of  canonically associated bundles (sometimes called {\it Tanaka's tower})
$$\ldots \longrightarrow P^{k+1} \overset{\pi_{k+1}}\longrightarrow P^{k} \overset{\pi_{k}}\longrightarrow P^{k-1}\overset{\pi_{k-1}}\longrightarrow \ldots \longrightarrow P^1 \overset{\pi_1} \longrightarrow P^0\overset{\pi} \longrightarrow M\ ,$$
where:
\begin{itemize}
\item[(i)] $\pi: P^0 \to M$ is the above introduced principal bundle with the structure group $G^0 = \Aut(\gm, J)$;
\item[(ii)] each bundle  $\pi_k:P^k\to P^{k-1}$ is a principal bundle with abelian structure  group $G^k$ with $Lie(G^k) \simeq \gg^k$
\end{itemize}
and for which the following holds:
 {\it If there exists
$k_o$ such that $\gg^{k_o +1} = \{0\}$,  then:
\begin{itemize}
\item[(1)] the CR automorphisms of $(M, \cD, J)$ form a Lie group of dimension $\dim \Aut(M, \cD, J) \leq \dim P^{k_o} = \dim(\gm + \gg^0 + \ldots + \gg^{k_o})$;
\item[(2)] there exists a canonical map which associates to each  $X \in \aut(M, \cD, J)$  a unique element  $\wh X^{k_o}$ in  a certain Lie algebra $\gp^{k_o}$ of vector fields of $P^{k_o}$, whose flows are local diffeomorphisms  preserving an appropriate  absolute parallelism on $P^{k_o}$;  the  correspondence $X \mapsto \wh X^{k_o}$  is a Lie algebra isomorphism  between $\aut(M, \cD, J)$ and $\gp^{k_o}$ and generalises  the above described   lifting map $X \to \wh X$, which transforms  infinitesimal CR automorphisms $X$ and vector fields on $P^0$;
\item[(3)] an element   $Y \in \gp^{k_o} (\simeq \aut(M, \cD, J))$ is identically vanishing    if and only if it vanishes at a single point     $y \in P^{k_o}$.
\end{itemize}
}
\par
If  there is $k_o\geq 0$ such that  $\gg^{k_o +1} = \{0\}$ holds, the CR structures of type $(\gm, J)$ are  called {\it of finite order} and the integer  $k_o$ is called the {\it order} of their   geometry.  \par
\smallskip
A (locally)  homogeneous CR manifold $(M, \cD, J)$ of type $(\gm, J)$ (\footnote{By ``locally homogeneous CR manifold''  we mean that the  Lie algebra $\aut(M, \cD, J)$  of infinitesimal CR transformations  generates local actions that are transitive  on open sets of the manifold.})  and with geometry of  finite order  $k_o$ is called a {\it maximally homogeneous  model} (resp. {\it maximally homogeneous local model}) of type $(\gm, J)$ if    $\Aut(M, \cD, J)$ (resp.  $\aut(M, \cD, J)$) has dimension equal to the dimension of $\gm + \gg^0 + \ldots + \gg^{k_o}$. \par
\smallskip
 By Prop. 10.7 in \cite{Ta},  if a fundamental CR algebra  $(\gm, J)$ corresponds to CR structures of  order $k_o$,  then  {\it the Lie algebra of infinitesimal  CR transformation  $\aut(G^\gm, \cD^\gm, J^\gm)$ of the standard CR model $G^\gm$  is isomorphic to $\gm + \gg^0 + \ldots + \gg^{k_o}$. In particular,  any  standard CR model  associated with structures   of finite order is a maximally homogeneous local model for  that class of geometric structures}.  Note that, under nondegeneracy conditions, a converse of this property have been proved in \cite{MN3}. We also remark that, as it was pointed out above,  each such  standard  model can be  realised as a regularly embedded submanifold of some $\bC^N$. \par
 \smallskip
The most familiar  examples of CR manifolds of finite order are probably  the real  hypersurfaces $M$ of $ \bC^N$ with  strictly positive  Levi forms:  their geometry is of order $k_o = 2$ and the unit sphere $S^{2N -1} = \SU_{N}/\U_{N-1} \times \bZ_2$ is a maximally homogenous model.  \par
\smallskip
We finally remark that, due to  the above properties (2) and (3) of Tanaka's towers,   if  the geometry of a   CR manifold has  order $k_o$, then  an infinitesimal CR transformation $X$ of such manifold is  zero if and only if  the corresponding lifted vector field  $X^{k_o}$  is equal to $0$  at some point $y \in P^{k_o}$ (hence, equal to $0$ at all points of $P^{k_o}$).
\medskip
\subsubsection{CR geometries of order $0$ and first order jets of automorphisms} In the very particular case of  a CR manifold $(M, \cD, J)$  with CR geometry of order $k_o = 0$,
the Lie algebra  isomorphism between $\aut(M, \cD, J)$  and the Lie algebra $\gp^0$ described in the above point  (2) of Tanaka's towers is nothing but  the canonical lifting map $X \to \wh X$ from infinitesimal CR automorphisms   and vector fields on $P^0$, as described in \S \ref{2.1.2}.  From this and the above  property (3),
it follows that  two infinitesimal CR automorphisms $X, X'$ coincide if and only if the lifted map $\wh{(X - X')}$ vanishes  at just one point $y \in P^0$ (and, hence, at all point of $P^0)$.  An explicit check shows that this occurs  if and only if there exists $x \in M$ such that the first order jet $j^1(X)|_x$ coincides with $j^1(X')|_x$.\par
\smallskip
Summing up,  we have that {\it for  manifolds with CR geometries of order $0$, the infinitesimal automorphisms are uniquely determined by their first order jets at a single point.}
\par
\smallskip
This property has immediate counterparts for  higher order CR geometries, as for instance  the   Levi nondegenerate real hypersurfaces.  For such different  kind of  CR manifolds,  the infinitesimal automorphisms are uniquely determined by jets of orders strictly  larger than one.
Due to this,  if  one looks for  analogies and/or differences between   $0$-th order CR geometries and  higher order ones,  the former  looks like  with  much fewer   degrees  of freedom  than  the  latter.  Such phenomenon is sometimes  mentioned   as the    {\it rigidity} feature of  $0$-th order geometries (see e.g. \cite{GK}).
\par
\medskip
\subsection{Complexifications of CR Tanaka algebras and  prolongations}
\subsubsection{Complexifications of Tanaka algebras}
Let $\gm$ be  a {\it real} fundamental  algebra  and  denote by $\gm^\bC = \gm + i \gm$  its associated complexified fundamental algebra.   Note   that the real subspace  $\gm \subset \gm^\bC$  can be characterised  as the fixed point set $(\gm^\bC)^\t$  of an appropriate anti-involution $\t$ of $\gm^\bC$, namely  of the conjugation    $\t(X + i Y) \=  \overline{X + i Y} = X - i Y$. This   yields   to the following   efficient characterisation of $\bC$-linear extensions of real linear  maps, whose proof is  basically  a straightforward  consequence of  the definitions:
\\[4pt]
(*)
 {\it A   map $L\in \Hom_\bC(\gm^\bC, \gm^\bC)$ is the $\bC$-linear extension of an $\bR$-linear map in  $\Hom_\bR(\gm, \gm)$ if and only if
$\overline{L(W)} = L(\overline W)$ for all $W = X + i Y \in \gm^\bC$.
}
\par
\smallskip
Consider now  a Tanaka algebra   $\gm + \gg^0$  over $\bR$. In this case the  natural map
$\imath: \Hom_\bR(\gm, \gm) \to$  $\Hom_{\bC}(\gm^\bC, \gm^\bC)$,   sending each $\bR$-linear map $L$  into its unique $\bC$-linear extension,  induces  a natural  injective  Lie algebras  homomorphism
$$\jmath^0 \= \imath|_{\gg^0}: \gg^0 \subset \aut(\gm) \longrightarrow \aut(\gm^\bC)\ . $$
 The complex subalgebra  $\Span_\bC(\jmath^0(\gg^0)) = \gg^0 + i \jmath^0(\gg^0)$ of $\aut(\gm^\bC)$ is clearly isomorphic to $(\gg^0)^\bC$, thus, for simplicity of notation, we  constantly identify  those two. Note that  the {\it complex} vector space $\gm^\bC + (\gg^0)^\bC$ has a natural structure of  {\it complex} Tanaka structure, associated with the    fundamental algebra  $\gm^\bC$,  and we call it {\it complexification} of $\gm + \gg^0$.\par
 \smallskip
The natural injections $\gm  \xhookrightarrow{\hskip 5pt}  \gm^\bC$,  $\gg^0  \xhookrightarrow{\hskip 5pt}  (\gg^0)^\bC$,  $\gg^i \xhookrightarrow{\hskip 5pt}  ((\gg^0)^\bC)^i$, $i \geq 1$ (the latter   determined by taking  $\bC$-linear extensions of the homomorphisms in the  $\gg^i$)  combine and determine an injective Lie algebras homomorphism
$$\jmath^\infty: (\gm + \gg^0)^\infty \longrightarrow (\gm^\bC + (\gg^0)^\bC)^\infty\ .$$
In this paper, we are particularly interested in  the injection
 $$\jmath^1 \= \jmath^\infty|_{\gg^1}: \gg^1 \subset (\gm + \gg^0)^\infty   \xhookrightarrow{\hskip 10pt}  \wh \gg^1 \= ((\gg^0)^\bC)^1 \subset (\gm^\bC +  (\gg^0)^\bC)^\infty\ .$$
The following technical lemma  gives a  characterisation of the image of $\jmath^1$.\par
\begin{lem}
\label{lemma31}
The subspace $\jmath^1(\gg^1)$ of  $\wh \gg^1$  consists of the $\bC$-linear maps  $L \in \Hom_\bC(\gm^\bC,  \gm^\bC + (\gg^0)^\bC)$ satisfying \eqref{defg1} for   $\ell = 1$ and   such that
\beq
\overline{L(W)(Z)} =  L(\overline W)(\overline Z)  \qquad \text{for all}\ \ W, Z \in \gm^\bC.
\eeq
\end{lem}
\begin{pfns} Note  that a map $L \in \wh \gg^1$ is the $\bC$-linear extension of  a linear map in  $\gg^1 \subset \Hom_\bR(\gm, \gm + \gg^0)$ if and only if it satisfies the following two conditions:
\begin{itemize}
\item[a)] The restriction $L|_{\gm}$ maps any element $X \in \gm^{-1}$ in $\jmath^0(\gg^0) \subset (\gg^0)^\bC$;
\item[b)] it satisfies the condition
\begin{multline*} [L(Y^{-1}), Z^{-p}] + [ Y^{-1}, L(Z^{-p})] = L([Y^{-1}, Z^{-p}])\\
\text{for any}\ \ Y^{-1} \in \gm^{-1}\ ,\ Z^{-p} \in \gm^{-p}\ ,\ -p = -\mu, \ldots, -1.
\end{multline*}
\end{itemize}
On the first hand, by definition of the space $\wh \gg^1$ and the fact that $\gm$ is generated by $\gm^{-1}$, we have that (a) implies (b), which can  be consequently neglected.  On the other hand, the elements of $\jmath^0(\gg^0)$ (i.e. the $\bC$-linear extensions of elements in $\gg^0$) are characterised by the above property (*).
Hence (a) holds if and only if
\beq \label{inter} \overline{L(X)(Z)} = L(X)(\overline Z)\qquad \text{for each}\ X \in \gm\ ,\  Z \in \gm^\bC\ .\eeq
This is equivalent to say that,  for all $W = X + i Y, Z \in \gm^\bC$,
$$ \overline{L(W)(Z)} = \overline{L(X)(Z)} + \overline{L(i Y)(Z)} =
\overline{L(X)(Z)}  - i  \overline{ L(Y)(Z)}  \overset{\eqref{inter}}= L(\overline W)(\overline Z)\ .\eqno\qed
$$
\end{pfns}
\par
\medskip
\subsubsection{Complexifications of   CR Tanaka algebras}
Let  $(\gm, J)$ be a  fundamental CR algebra and denote by  $\gm + \gg^0$,  $\gg^0 = \aut(\gm, J)$, the associated  CR Tanaka algebra.  Using the notation of previous subsection, we call {\it complexified fundamental CR algebra} the pair $(\gm^\bC, J)$ and we call
 {\it complexified CR Tanaka algebra of $(\gm, J)$} the complexified  Tanaka algebra  $\gm^\bC + (\gg^0)^\bC$. \par
 \smallskip
We recall that the $\bC$-linear extension  of the linear map  $J: \gm^{-1} \to \gm^{-1}$ on the complexification $(\gm^{-1})^\bC$
has   exactly two eigenvalues,  $+ i$ and  $-i$, and  two  associated eigenspaces,  $\gm^{10}$,  $\gm^{01} = \overline{\gm^{10}}$,  called {\it holomorphic} and {\it anti-holomorphic} subspaces. We also remind   that, conversely, for each direct sum decomposition $(\gm^{-1})^\bC  = \gm' + \gm''$ with $\gm'' = \overline{\gm'}$, there exists a unique complex structure $J$ on $\gm^{-1}$ for which the holomorphic and anti-holomorphic subspaces are precisely $\gm^{10} = \gm'$ and $\gm^{01} = \gm''$. Finally,  we observe  that the complex structure $J: \gm^{-1} \to \gm^{-1}$ satisfies the integrability condition \eqref{integrability} if and only if the corresponding holomorphic and anti-holomorphic subspaces $\gm^{10}, \gm^{01} \subset \gm^\bC$ satisfy the conditions
\beq  \label{integrability1}
[\gm^{10},  \gm^{10}]=  0 \qquad \text{and}\qquad [\gm^{01}, \gm^{01}] \left(= \overline{[\gm^{10}, \gm^{10}]} \right) =  0\ .
\eeq
Due to this, we have a natural one-to-one correspondence between fundamental CR algebras $(\gm, J)$ and pairs $(\gm^\bC, \gm^{10})$,  formed by
a complexified fundamental algebra $\gm^\bC$ and a subspace $\gm^{10} \subset (\gm^{-1})^\bC$, satisfying \eqref{integrability1} and such that the direct sum decomposition $(\gm^{-1})^\bC = \gm^{10}  + \overline{\gm^{10} }$ holds.\par
\smallskip
From this and  characterisation  (*)  of $\bC$-linear extensions, we get that  the $\bC$-linear extensions of the derivations $L \in \aut(\gm, J)$ are determined  by the following
\begin{lem}
\label{lem-C-linear-g0}
An element $L \in \aut(\gm^\bC)$ is  the $\bC$-linear extension of an element in $\gg^0 = \aut(\gm, J)$ if and only if $\overline{L(X)} = L(\overline{X})$ for any $X \in \gm^\bC$ and  $L(\gm^{10}) \subset \gm^{10}$.
\end{lem}
\par
\bigskip
\section{Totally nondegenerate CR manifolds}
\subsection{Universal fundamental CR algebras}
Given  a (real or complex) $n$-dimensional  vector space $V$ and an   integer $\mu \geq 1$, the following   result by Tanaka gives a unified presentation of
 all  fundamental  algebras $\gn$ of depth $\mu$ for which       $\gn^{-1} \simeq V$.
\begin{theo}\cite[\S 3]{Ta} \label{universal} For any finite-dimensional vector space $V$ over $\bK = \bR, \bC$ and any integer $\mu \geq 1$, there exists a   fundamental algebra $\gU_V$ of depth $\mu$, unique up to isomorphisms, satisfying the following two properties:
\begin{itemize}
\item[i)] $\gU_V^{-1} = V$;
\item[ii)] for any fundamental algebra $\gn$ of depth $\mu$  with $\gn^{-1} \simeq V$, there exists a {\rm natural} surjective Lie algebra homomorphism $\f_{(\gn)}: \gU_V \to \gn$, so that   $\gn   \simeq \gU_{V} /\gi$, with  $\gi = \ker \f_{(\gn)} $.
\end{itemize}
\end{theo}
In this statement,    the ``naturalness''  of the homomorphism  $\f_{(\gn)}$ has  the  following meaning: if $\psi: \gn \to \gn'$ is a Lie algebra homomorphism between two  fundamental algebras of same depth and same space of degree $-1$, and if  $\f_{(\gn)}: \gU_V \to \gn$ and $\f_{(\gn')}: \gU_V \to \gn'$ are the associated surjective homomorphisms, then
$\f_{(\gn')} = \psi \circ \f_{(\gn')}$  up to composition  with  isomorphisms.\par
\medskip
The graded Lie algebra $\gU_V$ of Theorem \ref{universal}  is called  {\it universal fundamental} (or {\it  free}) {\it Lie algebra of depth $\mu$  generated by $V$}. An explicit  construction of  $\gU_V$   is   given  in  Tanaka's proof of Theorem \ref{universal}  in  \cite[\S 3]{Ta} (see also \cite{Ha, MP, Wa}).
\par
\medskip
  Consider now a {\it real} vector space $V$ with a complex structure $J: V \to V$ and  let $V^\bC = V^{10}+ V^{01}$  be  the natural decomposition of $V^\bC$ as the direct sum of  $J$-eigenspaces. Further, given  $\mu \geq 1$, consider the  {\it complex} universal fundamental algebra  $ \gU_{V^\bC}$ of depth $\mu$ generated by $V^\bC$ and denote by $\gi_{10}, \gi_{01}$ the ideals of $\gU_{V^\bC}$, generated by the   $-2$-degree spaces $[V^{10}, V^{10}]$, $[V^{01}, V^{01}] $, respectively. Finally, let $\gU_{J}$ be  the quotient algebra
  \beq \label{universalCR1} \gU_{J} \= \gU_{V^\bC}/(\gi_{10} + \gi_{01})\ , \eeq
  which, by construction,  is a complex fundamental algebra.\par
  \smallskip
  Following the same  steps  of the proof in \cite{Ta}  of Theorem \ref{universal},  the following CR analogue to that result can be immediately derived.\par
  \begin{prop} \label{universalCR}  Given a $2n$-dimensional real vector space $V$  with  complex structure $J$  and an integer $\mu \geq 1$, the  algebra $\gU_{J}$ defined in \eqref{universalCR1} is a complexified   fundamental CR algebra with $((\gm^{-1})^\bC, J) = (V,^\bC J)$. \par
  Moreover,  for any other complexification $\gm'{}^\bC$  of a fundamental  CR  algebra $\gm'$ of depth $\mu$  with $((\gm'{}^{-1})^\bC, J')  \simeq (V^\bC, J)$, there is a {\rm natural} surjective Lie algebra homomorphism $\f_{(\gm^\bC)}: \gU_{J} \to \gm'{}^\bC$, so that
  $\gm'{}^\bC $ is isomorphic to the quotient   $\gU_{J} /\gi$  with   $\gi \= \ker \f_{(\gm^\bC)}$.
\end{prop}
Motivated by  this,  we call  $\gU_J$ the  {\it universal  fundamental CR algebra of depth $\mu$   generated  by $(V^\bC, J)$}.
 Note  that, since any complexified  fundamental CR algebra $\gm^\bC$ with $ (\gm^\bC)^{-1} = V^\bC$  is isomorphic to a quotient of  $\gU_J$ by an ideal $\gi \subset \sum_{i \geq 2} \gU_J^{-i}$, we have that $\gU_J$ can be characterised as {\it the  complexified  fundamental CR algebra, reaching the maximal values for the dimensions
 $n_k = \dim_\bC  \gU_J^{-k}$, $k \geq 2$,  amongst all complexified fundamental CR algebras having the same complex space  $V^\bC\= \gU^{-1}_J $ as  subspace of degree $-1$}.\par
 \medskip
 \subsection{Totally nondegenerate CR manifolds}\label{section32}
 We are now ready to give the definition of the Lie algebras, which are the  objects of study  of  this paper.  Let $\gm + \gg^0$ be a CR Tanaka algebra of depth $\mu$, with  complex structure $J: \gm^{-1} \to  \gm^{-1}$, and denote by $\gm^\bC + (\gg^0)^\bC$ the corresponding complexified CR Tanaka algebra. \par
\begin{definition}
\label{fundamentalB}  A (real) CR Tanaka  algebra  $\gm + \gg^0$ of depth $\mu \geq 3$  is called  {\it totally nondegenerate} if the  complexification  of $\gm$ has the form  $\gm^\bC = \gU_J/\gi^{-\mu}$ for an ideal    $\gi^{-\mu}$,  entirely  included in  the lowest degree subspace $\gU_J^{-\mu}$  of the universal  fundamental CR algebra $\gU_J$.\par
\end{definition}
\medskip
By the above characterisation of the  universal  fundamental CR algebras,   the totally nondegenerate CR Tanaka algebras can be also described as  the fundamental CR  algebras,  whose complexified  algebras $\gm^\bC$  reach the  maximal values for the dimensions   $n_k = \dim_\bC  (\gm^\bC)^{-k}$ for all   $k \geq 2$ {\it with the only possible exception of $k = \mu$}: in fact, the subspace $(\gm^\bC)^{-\mu}$ of lowest degree   might be a non-trivial quotient of $ \gU_J^{-k}$.\par
\smallskip
Note  that, by this property, if $\gm + \gg^0$ is totally nondegenerate, then  for each fixed $X^{-1}_o \in \gm^{-1}$, all linear maps
\beq \L^{(k)}: \gm^{-k} \rightarrow \gm^{-k-1}\ ,\quad  \L^{(k)}(Y^{-k}) \= [X^{-1}_o, Y^{-k}]\ ,\qquad  1 \leq k  \leq \mu-1\ ,
\eeq
 have  maximal possible rank or, equivalently, trivial kernels. This can be considered as a generalisation of  the classical  {\it Levi nondegeneracy condition} for  CR structures and motivates the name   ``totally nondegenerate'' for such fundamental CR algebras.
\par
\medskip
In correspondence  with  totally nondegenerate CR algebras we have the following
\begin{definition} \label{deftotallynon}  A   regular CR manifold  $(M, \cD, J)$ of type $(\gm, J)$ is called {\it totally nondegenerate} if the  associated Tanaka CR algebra $\gm + \gg^0$  is  totally nondegenerate.
\end{definition}
As we mentioned in the Introduction,  in \cite{Be1} Beloshapka introduced a special class of  embedded CR manifolds in  $\bC^N$, called {\it model surfaces}, each of them being totally nondegenerate.  He  proved that, for each model surface  $M$,   the    Lie algebra $\aut(M, \cD, J)$ of infinitesimal automorphisms surely includes  a subalgebra, which generates a locally transitive action on $M$  and  is isomorphic to  the   fundamental Lie algebra $\gm$  of such homogenous CR manifold. This means that the group $G^\gm$, which acts simply  transitively on the standard CR model (see \S \ref{2.1.2}), has a   simply transitive action on the (universal covering space of the)  model surface  with fundamental Lie algebra $\gm$. By a classical property of    Lie groups, for each fixed $o \in M$ and each linear map $\f: T_e G^\gm (\simeq \gm) \to T_o M$, there exists  a $G^\gm$-equivariant covering map $f^{(\f)}: G^\gm \to M$ such that $f^{(\f)}_*|_e = \f$. Hence, if $\f$ is a linear isomorphism between $\gm$ and $T_o M$ that gives also  an isomorphism between the  CR structures at $e \in G^\gm$ and $o \in M$,  respectively, using  the $G^\gm$-equivariance of $f^{(\f)}$, one can directly see  that   $f^{(\f)}$  maps the whole CR structure of $G^\gm$ onto the CR structure $M$. This means  that, {\it up to a covering,  $G^\gm$ and the corresponding Beloshapka's  model surface  $M$ are CR equivalent}.\par
\smallskip
 We conclude observing that, by Corollary 3 in \cite[\S 11]{Ta},  the geometry of totally nondegenerate CR manifolds is of finite order.  This immediately implies that for any Beloshapka's model $(M, \cD, J) $, the Lie algebra $\aut(M, \cD, J)$ of infinitesimal CR automorphisms is isomorphic to the  prolongation $(\gm + \gg^0)^\infty$ of  the Tanaka algebra $\gm + \gg^0$,  $\gg^0 = \aut(\gm, J)$ (\cite[Prop.10.7]{Ta}).  \par
 \medskip

\subsection{Statement of the main result}

\begin{theorem} \label{main} For any  totally nondegenerate CR  algebra $(\gm, J)$ of depth $\mu>3$, the space $\gg^1$ in the    prolongation $(\gm+\gg^0)^\infty$ of its Tanaka algebra is trivial.  Consequently, the CR geometry of any  totally nondegenerate CR manifolds  is of  order $0$.
\end{theorem}
As an immediate corollary  of this result,  by the remarks in \S \ref{section32},  we have that for any Beloshapka's model surface   of depth $\mu \geq 4$ (which we already observed to be  locally CR equivalent to a standard model $(G^\gm, \cD^\gm, J^\gm)$) the Lie algebra of infinitesimal CR automorphisms is isomorphic to $\gm + \gg^0$ with $\gg^0 = \aut(\gm, J)$. This proves the so-called ``maximum conjecture'' of \cite{Be2}  for all cases with  $\mu \geq 4$.
\par
\bigskip
\section{Linear independent sets in  totally nondegenerate  CR algebras}
Let $\gU_V$ be the universal fundamental algebra of depth $\mu$, generated by the $n$-dimensional vector space $V = \gU^{-1}_V$.
 For the sake of  simplicity, in this  section   we call {\it degree}  of a homogeneous element of $\gU_V$ the  {\it  opposite} of the integer which gives its degree as element of the negatively graded Lie $\gU_V$. So, for instance, in this section the elements in $\gU^{-1}_V$ are said to be of  {\it degree $+1$}, those in $\gU^{-2}_V$ are  of  {\it degree $+2$}, etc.
\subsection{Hall bases adapted to  universal fundamental CR algebras}
\setcounter{equation}{0}
\subsubsection{Hall bases of universal fundamental Lie algebras}
Given a (ordered)  basis $\cB = (E_1, \ldots, E_n)$ for  $V$,  we call   {\it  $\cB$-monomials of degree $1$} the  elements $E_i \in \cB \subset \gU_V^{-1}$ and, for $k \geq 2$,   we  call   {\it $\cB$-monomials of degree $k$} all  elements $X \in \gU_V^{-k}$
that are equal to  Lie brackets  $X = [Y, Z]$ between two  $\cB$-monomials of lower degree.   We may also say that  the  $\cB$-monomials
are the homogeneous elements of $\gU_V$ that can be  obtained  by taking  iterated  Lie brackets  between elements of the form $\ad_{E_{i_1}} \circ  \ad_{E_{i_2}} \circ \ldots \circ \ad_{E_{i_{r-1}}}(E_{i_r})$. This latter kind of   $\cB$-monomials   are called {\it elementary $\cB$-monomials}.
\par
\medskip
A  {\it Hall basis $\cH^\cB$  associated with   $\cB$}  is an (ordered) basis for $\gU_V$,  containing the $\cB$-monomials    selected as follows. For each degree $k\geq 1$, let us denote by $\cH^{\cB(k)}$ the subset of $\cH^\cB$ consisting of all its elements of degree $k$
and  denote by $\prec$ the total order relation  between such elements, which makes it an {\it ordered} basis. The sets $\cH^{\cB(k)}$ and the order relation are defined inductively by the following steps.\par
\begin{itemize}[itemsep=2pt, leftmargin=28pt]
\setlength\itemindent{10pt}
\item[{\bf  Step $\mathbf 1$.}] The subset $\cH^{\cB(1)}$ coincides with the set of all   $\cB$-monomials of degree $1$, i.e. all elements of $\cB \subset \gU^{-1}_V$.   The order relation $\prec$  between elements in  $\cH^{\cB(1)} $ is set  equal to the order of  $\cB$.
\setlength\itemindent{16pt}
\item[{\bf  Step $\mathbf{k_o}$.}]  Suppose that all sets  $\cH^{\cB (\ell)}$, $1 \leq \ell \leq k_o -1$,   have been defined and that the order   $\prec$  between  elements in the union $\bigcup_{1 \leq \ell \leq k_o-1} \cH^{\cB (\ell)}$  has been fixed in such a way  that if $\degree(X) < \degree(Y)$,   then $X \prec Y$.
Then the  set $\cH^{\cB(k_o)}$  is  defined as the collection of all $\cB$-monomials of degree $k_o$ of the form
$[U, V]$ with  $U, V \in  \cH^{\cB}$ and of degree lower  than $k_o$,   satisfying the following two conditions:
(i) $V \prec U$;  (ii) if  $U = [U_1, U_2]$ for some $U_i \in \cH^\cB$, then $U_2 \preceq V$.  The  ordering relation is also extended to  the larger set  $\bigcup_{1 \leq \ell \leq k_o} \cH^{\cB (\ell)}$  by setting: (a) for all $Y \in  \cH^{\cB(k_o)}$ and all $X$ of degree   less than $k_o$,  it is  set $X \prec Y$; (b)  if $X, Y$ are distinct, but  both in $\cH^{\cB(k_o)}$, then it is set either $X \prec Y$ or $X \succ Y$, according to an arbitrarily   chosen  order on the set of $\cB$-monomials of degree $k_o$.
\end{itemize}
\par
\smallskip
The ordered set $\cH^\cB =  \bigcup_{1 \leq \ell \leq \mu} \cH^{\cB (\ell)}$ is proven  to  be a basis for $\gU_V$. The proof of this, given  in \cite{Ha}, is based on the properties of an explicit procedure for each $\bX \in \gU_V$, which transforms  an expansion $\bX = \sum_{i=1}^N a_i X_i$  in terms of $\cB$-monomials $X_i$  into a  possibly new expansion $\bX = \sum_{\ell = 1}^{N'} a'_\ell X'_\ell$ in terms of $\cB$-monomials, all belonging to $\cH^\cB$. Such final expansion is proven to be independent from  the
original expansion  $\bX  = \sum_{i=1}^N a_i X_i$ of $\bX$.\par
\smallskip
For our purposes, it is convenient to shortly review such a procedure. It is defined  by induction on the degree $k$ of $\bX$. For $k = 1$, it simply consists in leaving the expansion $\bX = \sum_{i=1}^N a_i X_i$ unchanged. Indeed, all $\cB$-monomials of degree  $1$ are in $\cH^\cB$ and there is no need of changing the expansion.  Assume now that the procedure has been explicitly determined for each elements of degree $1 \leq k \leq k_o$ and assume  $\degree(\bX) = k_o +1$. In this case   the  final expansion  $\bX = \sum_{\ell = 1}^{N'} a'_\ell X'_\ell$  is reached by iteratively applying the following three steps: \par
\smallskip
\ \\
{\it First Step.}  Since each $\cB$-monomial $X_i$ has degree $k_o + 1 \geq 2$, it has  the form $[U_i, V_i]$ for some monomials $U_i, V_i$ of lower degree and for which the previously defined procedure allows to uniquely express them as  $U_i\=\sum_j a_{ij} U_j$ and $V_i\=\sum_k b_{ik} U_k$ for some elements    $U_j, U_k \in \cH^\cB$. In particular, we may re-write  $\bX$ in the form $\bX=\sum_{i,j,k} a_{ij} b_{ik} [U_j, U_k]$  with all $U_j, U_k \in \cH^\cB$.\\[10pt]
{\it Second Step .} In the expansion  $\bX=\sum_{i,j,k} a_{ij} b_{ik} [U_j, U_k]$, replace  each $\cB$-monomial  $[U_j, U_k]$  with  $U_j\prec U_k$ by the   expression $-[U_k, U_j]$.
\\[10pt]
{\it Third Step.} In the expansion $\bX=\sum_{i,j,k} a'_{ij} b'_{ik} [U_j, U_k]$ which is reached after  {\it Second Step}, consider  all terms
$[U_j, U_k]$ for which $U_j$ has the form $U_j:=[U'_m, U'_n]$ for some $U'_m, U'_n \in \cH^\cB$;  in case $U'_n \preceq U_k$, the $\cB$-monomial $[U_j, U_k]= [[U'_m, U'_n], U_k]$  is left unchanged, otherwise it is replaced by the expression  $[[U'_m, U_k], U'_n]-[[U'_n, U_k], U'_m]$.
\\[10pt]
As we already mentioned, by applying First Step, Second Step and Third  Step and then First Step, Second Step and Third  Step  again and so on, after  a finite number of iterations, one reaches an expansion $\bX = \sum_{\ell = 1}^{N'} X'_\ell$, in which all $X'_\ell$ are in $\cH^\cB$ and which is   stable under each of the above three steps (see \cite{Ha} for a detailed proof of this property).
\par
\subsubsection{Hall bases adapted to  a universal fundamental CR algebra} \label{compatible}
Let now $(V^\bC, J)$ be  the complexification of a real vector space $V$, with complex structure $J$, and  $\gU_J = \gU_{V^\bC}/\gi$, $\gi \= \gi_{10} + \gi_{01}$,  the universal CR fundamental algebra generated by $(V^\bC, J)$.  We recall that
 $\gi_{10}$ and  $\gi_{01}$ are  the ideals of $\gU_{V^\bC}$ generated by the  subspaces $[V^{10}, V^{10}]$ and $[V^{01}, V^{01}] $, respectively.
 \par
 \medskip
A basis $\cB = (E_i)_{1 \leq i \leq  2n}$ for $V^\bC = \gU^{-1}_{V^\bC}$ is called {\it adapted to the ideal $\gi_{10}$} (resp. {\it to the ideal $\gi_{01}$}) if the first $n$ elements $E_1, \ldots, E_n$ are  all in $V^{10}$ (resp. all in $V^{01}$). A   basis $\cB$ for $V^\bC$  is called {\it compatible  with $\gi$} if there are  two new orders on $\cB$,   the first    making $\cB$ adapted to $\gi_{10}$, the second  making it  adapted to $\gi_{01}$.\par
\begin{lem} \label{lemma41} If a basis $\cB  = (E_i)$  for $V^\bC = \gU^{-1}_{V^\bC}$  is adapted to $\gi_{10}$,  then for any associated Hall basis $\cH^\cB$ of $\gU_{V^\bC}$, one has that  $\cH^{\cB (10)} \= \cH^\cB \cap \gi_{10}$ is a  basis, as a vector space, for the ideal $ \gi_{10} $.  A similar property holds for bases adapted to $\gi_{01}$.
\end{lem}
\begin{pf}
 Let us denote by   $(A_j)_{j= 1}^N$   the (ordered)  set  of linearly independent elements  in  $\cH^{\cB (10)}  \= \cH^\cB \cap \gi_{10}$. We just need  to show that, for each $2 \leq k \leq \mu$, the subset  $\cH^{\cB (10|k)}  \= \cH^{\cB(10)} \cap \gU^{-k}_{V^\bC}$ of elements $A_j$ of degree $k$ is actually a basis  for the  intersection  $\gi^{-k}_{10} \= \gi_{10} \cap \gU_{V^\bC}^{-k}$.  We  prove this   by  induction on $k$. \par
  \smallskip
  For $k = 2$, the claim is true because the set $\cH^{\cB(10|2)}$  is not only linearly independent (it is a subset of a Hall basis) but it  also contains all   generators $ \{  [E_i, E_j]$, $1 \leq j  < i \leq n\} $ of $\gi_{10}$. Suppose now that the claim has been  proved for  $2 \leq k \leq k_o$ and let us show that  it  holds also for $k = k_o + 1$.
A non-trivial element $\bX \in \gi^{-(k_o +1)}_{10}$ has  the form $\bX = \sum_s[Y_s, Z_s]$
   for some $Y_s \in \gi_{10}$ and $Z_s \in \gU_{V^\bC}$.  The degrees of the $Y_s$ and $Z_s$ are   less than $k_o+1$.
  By the inductive hypothesis, each  $Y_s$ has the form
  $Y_s = \l^j_s A_j$
  for  monomials
  $A_j \in \bigcup_{\ell \geq 2}^{k_o} \cH^{\cB(10|k_o)}$, while each  $Z_s$ has the generic form
 $Z_s = \mu_s^r A_r + \nu_s^t B_t$ with $A_r \in \cH^{\cB(10)}$ and   $B_t \in \cH^\cB \setminus \cH^{\cB(10)}$.
 Thus
  \beq \label{theexpression}   \bX = \sum_s  \l_s^j \mu_s^r [A_j, A_r] +  \sum_s  \l_s^j \nu_s^t [A_j, B_t] \ .\eeq
The elements $[A_j, A_r]$,  $[A_j, B_t]$ appearing in \eqref{theexpression} are all  $\cB$-monomials and all in  $\gi_{10}$,  but not necessarily they are all in $\cH^\cB$. On the other hand,  we know that applying  a finite number of times   the sequence of the {\it First, Second and Third  Steps} of Hall's procedure  described above, one can  pass from \eqref{theexpression}   to  an equivalent   expansion, in which appear only $\cB$-monomials  in  $\cH^\cB$. We claim that  \\[10pt]
 {\it After each application of  the triple
``$\text{\it Second Step} \longrightarrow \text{\it Third Step}    \longrightarrow  \text{\it First Step}$''
 the  expansion of $\bX$  is still  a sum of $\cB$-monomials of the form $[A_i, A_j]$,  $[A_i, B_j]$, $[B_j, A_i]$.}\\[10pt]
To see this,  let us start with   a linear combination $\bX$ of $\cB$-monomials of the form $[A_i, A_j]$,  $[A_i, B_j]$, $[B_j, A_i]$ (as for instance, in \eqref{theexpression}).  The Second Step of Hall's algorithm  can only possibly exchange  $A_i$ with $A_j$ in  each monomial of the form $[A_i, A_j]$ or  exchange  $A_i$ with $B_j$ in all others. Hence,  after such step,  the form  of the expansion is  as claimed.  A more careful analysis is  needed for the Third Step, in which one has to replace  terms  of the form  $[[X, Y], Z]$,  for some
$[X, Y] = A_r$ and  $Z = A_m, B_n$,   or $[X, Y]   = B_s$ and $Z = A_m$,   in which  $Z \prec Y$.  We recall that  these terms must be replaced
by  the equivalent sum $- [[Y, Z], X] + [X, Z], Y]$.\par
\smallskip
We   want to show that,  whenever such  a replacement is needed,   {\it each  summand  in  $- [[Y, Z], X] + [X, Z], Y]$  is always  equal to  a  Lie bracket   $[W_1, W_2]$,  in which either $W_1$ or $W_2$ is in $\gi_{10}$}. \par
To see this,  consider at first the situations in which $Z = A_m \in \cH^{\cB(10)}$: in  such cases,   both  brackets $[Y,Z]$ and $[X, Z]$ are  in $\gi_{10}$ and the desired property  is  true. Second,  consider the  case in which $Z = B_m \in \cH^{\cB} \setminus \cH^{\cB(10)}$ and, consequently,   $[X, Y]$ is equal to some  $A_r \in  \cH^{\cB(10)} $ of  degree  $\wt k  \leq k_o$. Recall that  we are assuming that   $Z = B_m  \prec Y$.
In case $\wt k = 2$, the relation  $Z = B_m  \prec Y$ might occur only if $k_o + 1 = 3$ and all three elements $X, Y, B_m$ have degree $1$. Since we are assuming that $A_r = [X,Y]$ is in $\gi_{10} \cap \gU^{-2}_{V^\bC}$, this means that $X, Y$ are amongst the first $n$ elements of $\cB$, i.e.,  in $V^{10}$,  and $B_m \prec Y$ implies that also $B_m$ is in $V^{10}$. We thus have that also $[X, B_m]$, $[Y, B_m]$ are in $\gi_{10}$ and the desired claim holds. It remains to consider the cases $ 3 \leq \wt k \leq k_o$.  For those, we recall  that, by the iterative algorithm of construction of the monomials of the Hall basis,  the bracket  $[X,Y] = A_r$ is either a monomial  of
 the form $[X, Y]   = [A_m, A_n]$ or  a monomial of the form $[X, Y]   = [A_m, B_n]$  (or $[B_n, A_m]$).
In both  cases,  each summand of $- [[Y, Z], X] + [X, Z], Y]$ has the desired form: indeed, the first has the form $[W_1, W_2]$ with $W_2 = X = A_m \in \gi_{10}$, while the second has the form $[W_1, W_2]$ with $W_1 = [X, Z] = [A_m, Z] \in \gi_{10}$
(for the alternative cases in which $[X, Y]    = [B_n, A_m]$, one needs only to  exchange  the role of $X$ with the one of  $Y$).\par
\smallskip
Now that we know that,  after the Third Step,  we have a linear combination of brackets $[W_1, W_2]$ in which at least one $W_{i}$ is in $\gi_{10}$,  we  can newly perform  the  First Step (i.e.  express all elements $W_1$ and $W_2$  in terms of
the elements of the Hall Basis) and  obtain a new expansion of $X$,  which is once again a linear combination of elements of the form $[A_i, A_j]$,  $[A_i, B_j]$ or $[B_j, A_i]$
as claimed. This concludes the proof of the claim.\par
\medskip
As we  mentioned above, after a finite number of iterations of  these steps,  the expansion of  $\bX \in \gi_{10}^{-(k_o+1)}$
in terms of the elements of the Hall basis $\cH^\cB$ is reached.  But, by the above claim, such  final expansion consists of
elements of the form $[A_i, A_j]$,  $[A_i, B_j]$,  $[B_j, A_i]$, thus in $\gi_{10}$.  So,   $\gi^{-(k_o+1)}_{10}$ is generated by $\cH^{\cB(10|k_o +1)}$, which is therefore  a  basis for $\gi^{-(k_o+1)}_{10}$.
\end{pf}
\par
\medskip
\subsection{$\cB$-types and decompositions of universal CR fundamental algebras}
Let $\gU_{V^\bC}$ be a universal fundamental algebra of depth $\mu$ and $\cB = (E_i)_{i = 1, \ldots, 2n}$ a basis for $V^\bC = \gU_{V^\bC}^{-1}$. To  each $\cB$-monomial $X$,  we assign inductively an  element of  $\bN^{2n}$ as follows.  If $X$ is a $\cB$-monomial of degree $1$, i.e. a vector $X = E_i$ for some $ 1 \leq i \leq 2n$, we define {\it $\cB$-type of $X$} the $2n$-tuple
$$\cB\text{-}\type(E_i) = (0, \ldots, 0, \underset{i-\text{th entry}}1, 0, \ldots, 0)\ .$$
Suppose now that the type has been defined for all $\cB$-monomials of degree less than or equal to $k_o \geq 1$. Then if $X$ is a $\cB$-monomial of degree $k_o + 1$, thus of the form $X = [Y, Z]$ for  $\cB$-monomials of lower degree, we set
$$\cB\text{-}\type(X) = \cB\text{-}\type(Y) + \cB\text{-}\type(Z)\ .$$
Roughly speaking, the $\cB$-type is the $2n$-tuple, whose $i$-th entry gives   the number of times  the corresponding  vector $E_i$ is used in the construction of the $\cB$-monomial. Due to this property, if two bases $\cB$, $\cB'$ have the same elements, but differ by their order,  a  $\cB$-monomial $X$ is also a $\cB'$-monomial and  the $\cB$-type  and
 $\cB'$-type   differ only by an appropriate permutation of the entries.
\par
\smallskip
  Given  $I = (I_1, \ldots, I_{2n}) \in \bN^{2n}$, we denote   $|I| = \sum_{j = 1}^{2n} I_j$. In this way,  if $\cB\text{-}\type(X) = I$, then we have $\degree(X) =|I|$. \par
  \smallskip
  For any $I \in \bN^{2n}$ with  $0< |I| \leq \mu$, we  define
$$\gU^{-I}_{V^\bC}\= \Span_{\bC} \{\ {\rm monomials} \ X \in  \gU_{V^\bC}\ :\ \cB\text{-}\type(X) = I\ \}\ .$$
By the above described property of    types, two bases $\cB, \cB'$ of $V^\bC$,  which differ only by the order of their elements, determine the same subspaces $\gU^{-I}_{V^\bC}$, even if such space might be associated with  two   distinct $2n$-tuples $I$ and $I'$ of natural numbers (differing only by a permutation of the entries). Note also that,  for each $1 \leq k \leq \mu$, the vector subspace $\gU^{-k}_{V^\bC}$  admits  the direct sum decomposition
$$\gU^{-k}_{V^\bC} = \bigoplus_{|I| = k} \gU^{-I}_{V^\bC}\ .$$
\medskip
\begin{lem} Assume that the basis $\cB = (E_i)$ for $V^\bC$ is compatible with the ideal $\gi = \gi_{10} + \gi_{01}$ (see \S \ref{compatible} for the definition) and,  for each $I \in \bN^{2n}$ with $|I| \leq \mu$, let
$$\gi_{10}^{-I} \= \gi_{10} \cap \gU^{-I}_{V^\bC}\ ,\qquad \gi_{01}^{-I} \= \gi_{01} \cap \gU^{-I}_{V^\bC}\ .$$
 Then, for each $1 \leq k \leq \mu$ and each $2n$-tuple $I_o \in \bN^{2n}$, we have
\beq \label{42} \gi_{10} \cap \gU^{-k}_{V^\bC} = \bigoplus_{|I| = k} \gi^{-I}_{10}\ ,\qquad \gi_{01}  \cap \gU^{-k}_{V^\bC}  = \bigoplus_{|I| = k} \gi^{-I}_{01}\eeq
\beq \label{43} \gU^{-I_o}_{V^\bC}/ \gi \simeq  \gU^{-I_o}_{V^\bC}/\left(\gi_{10}^{-I_o} + \gi_{01}^{-I_o}\right)\ .\eeq
In particular,  the universal CR fundamental algebra $\gU_J = \gU_{V^\bC}/\gi$ admits the following  direct sum decomposition into vector subspaces
\beq \label{44} \gU_J \simeq  \bigoplus_{k = 1}^\mu\left( \bigoplus_{|I| = k} \gU^{-I}_{V^\bC}/ \left(\gi_{10}^{-I} + \gi_{01}^{-I}\right)\right)\ .\eeq
\end{lem}
\begin{pf} Let $\cB'$, $\cB''$ be two ordered bases for $V^\bC$, made with the elements in $\cB$ and ordered so that $\cB'$ is adapted to $\gi_{10}$ and $\cB''$ is adapted to $\gi_{01}$. Let also $\cH^{\cB'}$ and $\cH^{\cB''}$ be two associated Hall bases, with the corresponding subsets  $\cH^{\cB'(10)}$, $\cH^{\cB''(01)}$, which are  bases for $\gi^{10}$ and $\gi^{01}$ by Lemma \ref{lemma41}.   Finally, for each
$\cB'$\text{-}type $I'$, denote by $\cH^{\cB'(10)I'}$ the subcollections of $\cH^{\cB'(10)}$ consisting of $\cB'$-monomials of $\cB'$\text{-}type $I'$.  Similarly,    denote by $\cH^{\cB''(01)I''}$ the subcollections of $\cH^{\cB''(01)}$ consisting of $\cB''$-monomials of $\cB''$\text{-}type
$I''$.  Due to  the fact that  the basis $\cH^{\cB'(10)}$ of $\gi_{10}$  (resp. $\cH^{\cB''(01)}$ of $\gi_{01}$)    coincides with  the  union  of the subsets $\cH^{\cB'(10)I'}$ (resp. $\cH^{\cB''(01)I''}$), the direct sum decompositions \eqref{42} follow immediately.\par
 To prove \eqref{43}, let us denote by $I'_o$ and $I''_o$ the $\cB'$-type and $\cB''$-type, respectively, of  the monomials  which have $\cB$-type $I_o$. We recall that  $I'_o$ and $I''_o$ can be  derived from $I_o$ by  appropriate permutations of the entries and that
 $ \gU^{-I_o}_{V^\bC} =  \gU^{-I'_o}_{V^\bC} =  \gU^{-I''_o}_{V^\bC}$. So, using the direct sum decomposition \eqref{42}, we get  the  isomorphisms
$$\gU^{-I_o}_{V^\bC}/\gi_{10} \simeq \gU^{-I'_o}_{V^\bC}/\gi^{-I'_o} _{10} = \gU^{-I_o}_{V^\bC}/\gi^{-I_o} _{10} ,\qquad \gU^{-I_o}_{V^\bC}/\gi_{01} \simeq \gU^{-I''_o}_{V^\bC}/\gi^{-I''_o} _{01} =  \gU^{-I_o}_{V^\bC}/\gi^{-I_o} _{01}\ .$$
By standard  arguments, this  implies
\begin{multline*} \gU^{-I_o}_{V^\bC}/\gi = \gU^{-I_o}_{V^\bC}/(\gi_{10} + \gi_{01})   \simeq \\
\simeq  (\gU^{-I_o}_{V^\bC}/\gi_{10} )/ (\gi_{01}/\gi_{10})
 \simeq (\gU^{-I_o}_{V^\bC}/\gi^{-I_o}_{10} )/( \gi_{01}/\gi^{-I_o}_{10}) \simeq\\
 \simeq   \gU^{-I_o}_{V^\bC}/(\gi^{-I_o}_{10} + \gi_{01})
 \simeq (\gU^{-I_o}_{V^\bC}/ \gi_{01}) /(\gi^{-I_o}_{10}/\gi_{01})  \simeq\\
 \simeq  (\gU^{-I_o}_{V^\bC}/ \gi^{-I_o}_{01}) /(\gi^{-I_o}_{10}  /\gi^{-I_o}_{01})\simeq \gU^{-I_o}_{V^\bC}/(\gi^{-I_o}_{10} + \gi^{-I_o} _{01})\ ,
 \end{multline*}
  i.e. \eqref{43}.\par
\smallskip
The last claim is a direct consequence of \eqref{43}. In fact, the isomorphisms \eqref{43}  allow to construct an isomorphism between    each space $\gU^{-k}_{V^\bC}/\gi \subset \gU_J$  and the direct sum of all quotients $\gU^{-I}_{V^\bC}/(\gi^{-I} _{10} + \gi^{-I} _{01})$ with $|I| = k$.
\end{pf}
This lemma  yields to the following
\begin{cor}  \label{cor43} Let $\cB$ be a basis for $V^\bC$, which is compatible with $\gi = \gi_{10} + \gi_{01}$.  If
 $X_1, \ldots, X_\ell \in \gU_{V^\bC}$  are    $\cB$-monomials, whose $\cB$-types are  pairwise distinct and none of them is  in $\gi$, the  classes  $\wh X_i = X_i + \gi$  in $\gU_J = \gU_{V^\bC}/\gi$   constitute a   linearly independent  set  in $\gU_J$.
\end{cor}
\begin{pf} For each $i$, let $I_i$ be the $\cB$-type of $X_i$. Since the types $I_i$ are pairwise distinct,   the classes $\wh X_i$  belong to the distinct subspaces $\gU^{-I_i}_{V^\bC}/ \left(\gi_{10}^{-I_i} + \gi_{01}^{-I_i}\right)$ and are therefore linearly independent, because of   the direct sum decomposition \eqref{44}.
\end{pf}
\par
\medskip
\subsection{Linear independent sets in a totally nondegenerate  CR algebra}
We recall that a  (real) CR Tanaka  algebra
$$\gm + \gg^0 = \gm^{-\mu} + \gm^{-\mu + 1} + \ldots + \gm^{-1} + \gg^0$$
with complex structure $J$ on $\gm^{-1}$  is called   {\it totally nondegenerate} if and only if   $\gm^\bC = \gU_J/\gi^{-\mu}$ for some ideal    $\gi^{-\mu}$, which is entirely  included in the subspace $\gU_J^{-\mu}$  of the universal  fundamental CR algebra $\gU_J$ of depth $\mu$.
 This means that each homogeneous element of $\gm^\bC$
 of degree  strictly less than $\mu$ can be  identified with an   element of the universal algebra $\gU_J$, while the elements in  $(\gm^{-\mu})^\bC$ are  identifiable with  equivalence classes of  $\gU^{-\mu}_J/\gi^{-\mu}$. \par
 On the other hand, the elements of the universal CR fundamental algebra $\gU_J$ are in turn  classes of the quotient $\gU_{V^\bC}/(\gi_{10} + \gi_{01})$. \par
 \smallskip
 Since many  of the following arguments are based on  properties of the elements of the  algebra $\gU_J = \gU_{V^\bC}/(\gi_{10} + \gi_{01})$ and on the corresponding elements in $\gm^\bC   = \gU_J/\gi^{-\mu}$, for the sake of clarity, it is crucial  to fix a   notational convention,  which allows to easily indicate  if an element is  in $\gU_J$ or $\gm^\bC = \gU_J/\gi^{-\mu}$. So, from now on,
 \begin{itemize}[itemsep=2pt, leftmargin=28pt]
\item[--] the elements of $\gU_{J}$ are     denoted by capital roman letters as $X, Y, Z, \ldots$,
\item[--] the corresponding  images in $\gm^\bC = \gU_J/\gi^{-\mu}$ are indicated by $X^{\gm}, Y^{\gm}$,  etc.
\end{itemize}
 By previous remarks, a homogeneous element $X \in \gU_J$  can be different from the corresponding $X^\gm \in \gm^\bC$ if and only if $X \in \gU^{-\mu}_J$.  Due to this,  the  elements  of $\sum_{k = 1}^{\mu - 1} (\gm^{-k})^\bC$ will be often tacitly  identified with the corresponding elements in $\gU_J$ and hence denoted by capital roman letters as $X, Y, Z, \ldots$. The above   special notation $X^\gm$, $Y^\gm$, $\ldots$, for elements in $\gm^\bC$ is going to be used very rarely, namely only when  elements in $(\gm^{-\mu})^\bC$ are considered  and it is important to prevent misunderstandings.\par
 \medskip

 Corollary \ref{cor43} has the following immediate consequence
\begin{cor}  \label{cor44} Let $\cB = (E_i)$ be a basis for $V^\bC = (\gm^{-1})^\bC = \gU^{-1}_{V^\bC}$, which is compatible with $\gi = \gi_{10} + \gi_{01}$.  If
 $\wt X_1, \ldots, \wt X_\ell \in \gU_{V^\bC}$  are   non-zero $\cB$-monomials of degree strictly less than $\mu$,  with   pairwise distinct   $\cB$-types and none of them in $\gi$, the corresponding projected  elements   $X_i ( = X^{\tiny \gm}_i) \in \gm^\bC = \gU_{V^\bC}/(\gi_{10} + \gi_{01} + \gi^{-\mu})=\gU_J/\gi^{-\mu}$   constitute a   linearly independent  set  in $\gm^\bC$.
\end{cor}
\par
\bigskip
\section{Proof of  Theorem \ref{main}}
\subsection{Notational issues}
Let $\gm^\bC + (\gg^0)^\bC = \gU_J/\gi^{-\mu} + (\gg^0)^\bC$ be the complexification of a totally nondegenerate CR algebra. As  above,
we denote by  $X, Y,  \ldots$ elements in $\gU_J$ and we tacitly identify them with the corresponding elements in $\gm^\bC$ whenever
 they belong to  $\sum_{k = 1}^{\mu - 1} \gU^{-k}_J$. Otherwise, their equivalence classes in $\gm^\bC$  are denoted by $X^{\gm}$, $Y^{\gm}$, etc.
We point out that  any element  in $\gg^1$ will be also constantly and tacitly  identified with the corresponding element,  obtained by $\bC$-linearity extension,  in the subspace $\jmath(\gg^1) \subset \wh \gg^1$ of the prolongation of the complexification $\gm^\bC + (\gg^0)^\bC$.
\par
\medskip
Given $X, Y \in \gm^\bC$ and   $L \in  \gg^1 (= \jmath(\gg^1))$, we set:
\begin{itemize}
\item[--]  $X^r\bdot Y \= (\ad_X)^r(Y)$ for each integer $r \geq 0$ (we assume   $(\ad_X)^0 (Y) \=  Y$);
\item[--] $L_X \=  [L, X]$;
\item[--]  $L_{X|Y} \= [[L, X], Y] $.
\end{itemize}
According to such  convention,   for any  $L \in \gg^1$,   $E \in (\gm^{-1})^\bC$ and  $W \in ( \gm^{-k})^\bC$,
$$L_E \in (\gg^0)^\bC\ ,\qquad L_{E|W} \in (\gm^{-k})^\bC\ , \qquad L_{E^r \bdot W} \in (\gm^{-k - r + 1})^\bC\ .$$
\par \medskip
 \subsection{A basic property of the elements  $L \in \gg^1$}
Recall that,  by the integrability condition on $J$ and  Lemma \ref{lem-C-linear-g0}, for each  $L \in \gg^1$
 and  for any   $E, F \in \gm^{10}$ (so that their conjugate elements $\overline E, \overline F$ are in $\gm^{01}$)
  \begin{align}
&\label{5.5} [E, F] = 0\qquad \text{and} \qquad  [\overline E, \overline F] = 0 \ ,\\
\label{5.6} & L_{E|F} \in \gm^{10}\qquad \text{and}\qquad L_{\overline E|F} \in \gm^{10}\ ,\\
& L_{E|\overline F}  = \overline{L_{\overline E|F}} \in \gm^{01}\qquad \text{and}\qquad L_{\overline E|\overline F} = \overline{L_{E|F}} \in \gm^{01}\ .
\end{align}
Using these basic identities, we get the next useful lemma, which  can be considered as a CR version of \cite[Lemma 1]{Wa}.
\begin{lem}
\label{cor-Warhurst1}
Given  $E \in \gm^{10}$ and $L \in \gg^1$,  for any other  element $W \in \gm^\bC$
\begin{align} \label{war1-CR}  &L_{E\bdot W} =  L_{E|W} -  L_{W|E}  \ ,
\end{align}
and, more generally, for $r \geq 2$,
\beq  \label{war1-CRbis}
L_{E^r \bdot W} =  E^{r-1} \bdot  \left( r L_{E|W}  -   L_{W|E}  \right)  +  \frac{r (r-1)}{2}
E^{r-2}\bdot (L_{E|E} \bdot   W)\ .
\eeq
\end{lem}
\begin{pf} We first observe that, from \eqref{5.6} and \eqref{5.5}
\beq \label{comm} [\ad_{L_{E|E}},  \ad_E] = \ad_{[L_{E|E}, E]} = 0\ ,\eeq
meaning that the adjoint operators $\ad_{L_{E|E}}$ and $\ad_E$  commute.
Using this, we claim that for any  $\ell \geq 2$,
\beq\label{intermidiate} [L_E,   E^{\ell} \cdot W] = \ell E^{\ell-1} \bdot (L_{E|E} \bdot W ) + E^\ell \bdot L_{E|W}\ .\eeq
Indeed, when $\ell = 2$ this is true because
\begin{multline*} [L_E, [E, [E, W]]] \overset{\text{Jacobi id.}}= L_{E|E} \bdot( E \bdot W) + E \bdot ([L_E, [E, W]]) \overset{\eqref{comm}\ \&\ \text{Jacobi id.}}= \\
 = E \bdot (  L_{E|E} \bdot W)  +  E \bdot (L_{E|E} \bdot W) +
 E^2 \bdot L_{E|W}\ .\end{multline*}
 Assume now that \eqref{intermidiate} has been proved  for  $2\leq \ell \leq \ell_o-1$.  Then by \eqref{comm}, since   $\ad_{L_{E|E}}$ and $\ad_E$ commute, by  the inductive hypothesis,
\begin{multline*}[L_E, E^{\ell_o} \bdot W] \overset{\text{Jacobi id.}}= L_{E|E} \bdot (E^{\ell_o-1} \bdot W) + E \bdot [L_E,  E^{\ell_o-1} \cdot W] = \\
 = E^{\ell_o-1} \bdot (L_{E|E} \bdot W) +  (\ell_o -1) E^{\ell_o-1} \bdot (L_{E|E} \bdot W ) + E^{\ell_o} \bdot L_{E|W}=\\
  = \ell_o E^{\ell_o-1} \bdot (L_{E|E} \bdot W ) + E^{\ell_o} \bdot L_{E|W}\ .  \end{multline*}
We can now prove the lemma. We first observe that  \eqref{war1-CR} is nothing but the Jacobi identity.  For  \eqref{war1-CRbis}, the case $r = 2$ is obtained
by  iterated uses of  Jacobi identity. Indeed,
\begin{multline*}
L_{E^2 \bdot W} = [L, [E, [E, W]]] = [[L, E], [E, W]] + [E, [L, [E, W]]] =\\
=[ [[L, E],E], W] + [E, [[L, E],W]] +  E\bdot L_{E \bdot W}\overset{\eqref{war1-CR}}= \\
= 2   E\bdot  L_{E|W} -  E \bdot L_{W|E} + L_{E|E} \bdot W \ .
\end{multline*}
Assume now that  \eqref{war1-CRbis}  holds  for $2 \leq r \leq r_o-1$. Then,  by  \eqref{intermidiate}  and induction
\begin{multline*}
L_{E^{r_o} \bdot W} = [L_E , E^{r_o-1}\bdot W]+ E \bdot (L_{E^{r_o-1} \bdot W}) =\\
= (r_o-1) E^{r_o-2} \bdot (L_{E|E} \bdot W ) + E^{r_o-1} \bdot L_{E|W}
+ E^{r_o-1} \bdot  \left( (r_o-1) L_{E|W}  -   L_{W|E}  \right)  + \\
+ \frac{(r_o-1)(r_o -2)}{2}
 E^{r_o-2}\bdot (L_{E|E} \bdot  W)\ ,
\end{multline*}
 which  gives \eqref{war1-CRbis} for $r = r_o$.
\end{pf}\par
\medskip
\subsection{A  finer analysis of the properties of the elements   $L \in \gg^1$.}
\begin{lem} \label{firstlemma} If $\mu \geq 4$, then  for each $L \in \gg^1$ and  each  $E \in \gm^{10}$,  there exists  $\l^{(E)} \in \bC$ such that
$L_{E|E} = \l^{(E)} E$.
\end{lem}
\begin{pf} Without loss of generality, we can assume that $E\neq 0$. Let $\r$ be an eigenvalue
of the $\bC$-linear map  $L_E|_{\gm^{01}}: \gm^{01} \to \gm^{01}$ and  $\overline F_o$  a $\r$-eigenvector.
Since $E^{\mu} \cdot \overline F_o$ has degree $\mu+1$ and the depth of $\gm$ is $\mu$, we have that $E^{\mu} \cdot \overline F_o = 0$. Consequently,  setting $\wh E \=  L_{E|E}$,  from  \eqref{5.5}, \eqref{5.6},   \eqref{war1-CRbis} and the hypothesis   $L_{E|\overline F_o} = \r \overline F_o$,
\beq\label{5.12} 0 = L_{E^{\mu} \cdot \overline F_o} =  \r \mu   E^{\mu -1} \bdot  \overline{F}_o + \frac{\mu (\mu -1)}{2} E^{\mu -2} \bdot (\wh E\bdot \overline{F}_o)\ .\eeq
Applying $L$ to both sides  and,  once again,  using     \eqref{5.5}, \eqref{5.6},   \eqref{war1-CRbis} and   $L_{E|\overline F_o} = \r \overline F_o$ we get
\begin{multline} \label{3.6-1}0 = L_{(L_{E^\mu \bdot \overline F_o})} = \\
=  \r^2\mu (\mu -1) E^{\mu -2} \bdot  \overline{F}_o + \r \frac{\mu (\mu - 1) (\mu -2)}{2} E^{\mu -3} \bdot(\wh E \bdot \overline F_o)  +\\
+ \frac{\mu (\mu -1) (\mu-2)}{2} E^{\mu -3}  \bdot L_{E| ( \wh E \bdot \overline F_o)} -  \frac{\mu (\mu -1) }{2} E^{\mu -3}  \bdot L_{ ( \wh E \bdot \overline F_o)| E} +  \\
 +   \frac{\mu (\mu -1)(\mu - 2) (\mu -3)}{4}
E^{\mu -4}   \bdot ( \wh E \bdot(\wh E \bdot \overline F_o))   =\\
=  \r^2\mu (\mu -1) E^{\mu -2} \bdot  \overline{F}_o + \r \frac{\mu (\mu - 1) (\mu -2)}{2} E^{\mu -3} \bdot(\wh E \bdot \overline F_o)  +\\
+ \frac{\mu (\mu -1) (\mu-2)}{2} E^{\mu -3}  \bdot( L_{E|\wh E} \bdot \overline F_o) + \ \frac{\mu (\mu -1) (\mu-2)}{2} E^{\mu -3}  \bdot ( \wh E \bdot L_{E|\overline F_o} )+ \\
+ \frac{\mu (\mu -1) }{2} E^{\mu -2}  \bdot L_{ \wh E \bdot \overline F_o} +  \\
 +   \frac{\mu (\mu -1)(\mu - 2) (\mu -3)}{4}
E^{\mu -4}   \bdot ( \wh E \bdot(\wh E \bdot \overline F_o))
 \ .\end{multline}
Suppose now that $\wh E = L_{E|E}$ is not   a multiple  of  $E$. We may therefore
 consider a basis $\cB = (E_i, \overline F_j)$ for $\gm^{-1}$, in which:  (a) $(E_i)$ is a basis for  $\gm^{10}$  with $E_1 \= E$,  $E_2 = \wh E$, (b)  $(\overline F_j)$ is a basis for  $\gm^{01}$  with $F_1 \= F_o$.  Identifying  the elements of this  basis  with the corresponding elements in $\gU^{-1}_{V^\bC} \simeq (\gm^{-1})^\bC$, we have  that  $\cB$ is   compatible with the ideal $\gi = \gi_{10} + \gi_{01}$ (see \S \ref{compatible}).
Moreover, since $\mu \geq 4$, the last summand in the right hand side of \eqref{3.6-1} is a non-trivial $\cB$-monomial with  $\cB$-type
$$(\mu -4, 2 , 0, \ldots, 0, 1, 0, \ldots, 0)\ .$$
On the other hand, expanding all other summands in \eqref{3.6-1} in terms of $\cB$-monomials, we see that all terms in such  expansion
are $\cB$-monomials of $\cB$-type of the form
$$ (\mu - 2, *,*   \ldots, *)\qquad \text{or}\qquad (\mu - 3, *,*   \ldots, *)\ .$$
Hence, by  Corollary \ref{cor44}, the last summand cannot be equal to  a  linear combination   of the other summands and this contradicts   \eqref{3.6-1}. \end{pf}

\begin{rem} \label{firstlemmabis} Note that the statement of Lemma \ref{firstlemma} holds also when $\mu = 3$ and either  $\dim_\bC \gm^{10}  = 1$ or $\gm^\bC = \gU_J$. Indeed, if $\dim_\bC \gm^{10}  = 1$, then $L_{E|E}$ is proportional to $E$ simply because $\gm^{10}$
admits no more than one linearly independent vector. In case $\gm^\bC = \gU_J$, that is if $\gi^{-\mu} = \{0\}$,  then \eqref{5.12} is  satisfied in $\gm^\bC$ if and only if it is satisfied in $\gU_J^{-\mu}$. By the same arguments of the lemma, based on comparison of $\cB$-types of the summands, \eqref{5.12} is satisfied in $\gU_J$ if and only if $\wh E = L_{E|E} = - \frac{\r(\mu-1)}{2} E$.
\end{rem}
\par \medskip
In all the following,  {\it we constantly and tacitly assume that $\gm + \gg^0$ is a totally nondegenerate CR Tanaka algebra of depth $\mu \geq 4$}, so that Lemma \ref{firstlemma} surely applies to such Tanaka algebra. However, note that,  due to Remark \ref{firstlemmabis},  all results of this section are valid also in case the  depth is $\mu = 3$ and  either  $\dim \gm^{10} = 1$ or  $\gm^\bC = \gU_J$. \par
\medskip

\begin{lem}\label{secondlemma}  For any   $L \in \gg^1$,  the map
$$\a_L: \gm^{10} \to \bC\ ,\qquad \a_L(E) \=  \l^{(E)}\ ,$$
with $\l^{(E)}$ as   in Lemma \ref{firstlemma},  is  linear. Moreover,
for any pair $E, E' \in \gm^{10}$
\beq \label{5.13} L_{E|E'} = \frac{1}{2} \a_L(E') E +  \frac{1}{2} \a_L(E) E'\ .\eeq
\end{lem}
\begin{pf} Given $a \in \bC$ and $E \in \gm^{10}$,
$$\l^{(aE)} a E = L_{a E| aE} = a^2 L_{E|E} = a^2 \l^{(E)} E\ .$$
This means that $\l^{(aE)} = a \l^{(E)}$ for each $a \in \bC$. Thus, in order to prove that $\l^{(E)}$ is linear in $E$, it remains to prove that
$\l^{(E + E')} = \l^{(E)} + \l^{(E')}$ for any pair of linearly independent vectors $E, E' \in \gm^{10}$. For this, we first note that, by \eqref{war1-CR} and \eqref{5.5}
$0 = L_{E \bdot E'} = L_{E|E'} - L_{E'|E}$, meaning that $L_{\cdot|\cdot}$ is symmetric in its arguments if both of them are in $\gm^{10}$. Moreover,
$$\l^{(E+E')} (E + E') = L_{E+ E'| E + E'} = \l^{(E)} E + 2 L_{E|E'} + \l^{(E')} E'\ ,$$
which implies that
\beq \label{ecco} L_{E|E'} = \frac{1}{2} (\l^{(E+E')} - \l^{(E)}) E + \frac{1}{2} (\l^{(E+E')} - \l^{(E')}) E'\ .\eeq
Now, by bilinearity of $L_{\cdot|\cdot}$, we observe that
\begin{multline*}
\frac{1}{2} (\l^{(2 E+E')} - \l^{(2E)}) 2E + \frac{1}{2} (\l^{(2E+E')} - \l^{(E')}) E' = \\
= L_{2E|E'} = L_{E|2 E'}  =\\
=  \frac{1}{2} (\l^{(E+ 2 E')} - \l^{(E)}) E + \frac{1}{2} (\l^{(E+2 E')} - \l^{(2 E')}) 2 E'\ .
\end{multline*}
Since $E$ and $E'$ are linearly independent and $\l^{(2 E)} = 2 \l^{(E)}$, we infer that
\beq
\nonumber
\begin{split} & 2 \l^{(2 E+E')} -  3 \l^{(E)}  - \l^{(E+ 2 E')}   = 0\ ,\\
&\l^{(2E+E')}  - 2 \l^{(E+2 E')} + 3 \l^{(E')} = 0\ .
\end{split}
\eeq
Subtracting twice the second  from the first equation, we get
$ -  3 \l^{(E)}  + 3 \l^{(E+2 E')} - 6 \l^{(E')} = 0$, i.e.  $\l^{(E + 2 E')} = \l^{(E)} + \l^{(2 E')}$,  which  proves the linearity of $\a_L$. The last claim is a direct consequence of \eqref{ecco} and linearity of $\a_L$.
\end{pf}

\begin{lem} \label{32} Let $L \in \gg^1$ with associated linear map $\a = \a_L \in (\gm^{10})^*$. Then, for each $0 \neq E \in \gm^{10}$ and $0 \neq \overline F \in \gm^{01}$,
\beq \label{res}  (\ad_{L_E})^2 (\overline F)  = -  \a(E) \frac{2 \mu -3}{2}  \ad_{L_E}(\overline F)  -   \a(E)^2  \frac{(\mu -1)(\mu - 2) }{4}    \overline F. \eeq
\end{lem}
\begin{pf}  Given $\overline F$, let $\nu$ be  the first integer such that $E^\nu \bdot \overline F = 0$. By definition of totally nondegenerate CR algebras,    $\nu$  is  either $\mu-1$ or $\mu$.   We also  set   $\overline F' \= L_{E|\overline F} = \ad_{L_E}(\overline F)$ and, as usual,  we identify
$E$, $\overline F$,  $\overline F'$ with the corresponding elements in $\gU^{-1}_{V^\bC} = (\gm^{-1})^\bC = V^\bC$. Finally, we denote by $\cB = (E_i, \overline F_j)$ a  basis for $\gU_{V^\bC}^{-1}$, compatible with $\gi = \gi_{10} + \gi_{01}$, with $E_1 \= E$ and $\overline F_1 \= \overline F$. In case $\overline F'$ is not a multiple of $\overline F$, we  also assume that $\overline F_2 = \overline F'$.\par
\medskip
We now prove \eqref{res} by considering the two possibilities for $\nu$. \\[10pt]
\noindent{\it Case 1: \ $\nu = \mu -1$}. \\[5pt]
 Since $E^{\mu -1} \bdot \overline F = 0$,  from \eqref{5.5}, \eqref{5.6}, \eqref{war1-CRbis} and  \eqref{5.13},
 $$0= L_{E^{\mu-1} \bdot \overline{F}}  = (\mu -1) E^{\mu -2} \bdot  \overline{F}'
+ \a(E) \frac{(\mu - 1) (\mu -2)}{2} E^{\mu -2}  \bdot \overline{F} \ .
$$
Since all terms of this equality have degree $\mu -1$, they can be identified with corresponding elements in $\gU_J$. Considering the $\cB$-type of the two $\cB$-monomials involved, by  Corollary \ref{cor44}  we have that the equality holds if and only if
$$\ad_{L_E}(\overline F) =  \overline F' = -  \a(E)\frac{\mu -2}{2}  \overline F\ .$$
Replacing this into  \eqref{res}, we see that the equality is satisfied, proving the lemma in this case.  \par
\medskip
\noindent{\it Case 2: \  $\nu = \mu$}. \\
Since $E^{\mu} \cdot \overline F = 0$, using as usual \eqref{5.5}, \eqref{5.6} and  \eqref{war1-CRbis}, we have
$$ 0 = L_{E^{\mu} \cdot \overline F} =\mu  E^{\mu -1} \bdot  \overline{F}'  + \a(E) \frac{\mu (\mu -1)}{2} E^{\mu -1} \bdot \overline{ F} \ .$$
 Dividing by $\mu$,  applying $L$ to both sides and  using once again \eqref{5.5}, \eqref{5.6} and  \eqref{war1-CRbis}, we get
 \begin{multline*}
0 = (\mu-1) E^{\mu-2} \bdot L_{E|\overline F'} + \a(E) \frac{(\mu-1) (\mu -2)}{2} E^{\mu-2} \bdot \overline F' + \\
+ \a(E) \frac{(\mu-1)^2}{2}E^{\mu-2} \bdot \overline F' + \a(E)^2 \frac{(\mu-1)^2(\mu -2)}{4} E^{\mu-2} \bdot \overline F \ .
\end{multline*}
With the usual arguments based on $\cB$-types, we get  that this can be satisfied if and only if  \eqref{res} holds.
\end{pf}
\begin{prop}\label{cor1} Let $L \in \gg^1$ with associated linear map $\a = \a_L \in (\gm^{10})^*$.
Then, for each $E \in \gm^{10}$, the restricted map $\ad_{L_E}|_{\gm^{01}}: \gm^{01} \to \gm^{01}$ is diagonalisable, with at most  two eigenvalues, whose values  can be only
\beq \label{possible-eigen}\r = - \a(E) \frac{\mu-1}{2}\qquad \text{or}\qquad \r =  - \a(E) \frac{\mu-2}{2}\ .\eeq
In particular,  if $E \in \ker \a$, then $\ad_{L_E}|_{\gm^{01}} = 0$.
\end{prop}
\begin{pf} Let $\r$ be an eigenvalue of  $\ad_{L_E}|_{\gm^{01}}: \gm^{01} \to \gm^{01}$  and  $\overline F_o \neq 0$ an associated eigenvector. Then, by \eqref{res},
\begin{multline*} 0 =  \left(\r^2  +   \a(E) \frac{2 \mu -3}{2}  \r +  \a(E)^2  \frac{(\mu -1)(\mu - 2) }{4} \right) \overline F_o =\\
=  \left(\r + \a(E) \frac{\mu -1}{2}\right) \left(\r +  \a(E) \frac{\mu -2}{2}\right) \overline F_o\ .\end{multline*}
This implies  that \eqref{possible-eigen} are the only possibilities for the eigenvalues.
\par
\medskip
We now claim that only the following two  cases
might occur:
\begin{itemize}
\item[a)] $\a(E) \neq 0$ and  the linear map  $\ad_{L_E}|_{\gm^{01}}$ admits a basis for $\gm^{01}$ made of eigenvectors;
\item[b)] $\a(E) = 0$ and hence with all eigenvalues of    $\ad_{L_E}|_{\gm^{01}}$     equal to $0$.
\end{itemize}
For checking this, assume that neither  (a) nor (b)  occurs, i.e. that $\a(E) \neq 0$ and that $\ad_{L_E}|_{\gm^{01}}$ is not diagonalisable. This means that there is a basis $(\overline F_1, \ldots, \overline F_n)$ for $\gm^{01}$,  in which the associated matrix of $\ad_{L_E}|_{\gm^{01}}$ is in Jordan canonical form with at least one Jordan  block of order greater than or equal to $2$.  However, by \eqref{res}, no Jordan block can be of order strictly larger than $2$. We may therefore assume that the basis $(\overline F_j)$ is such that
$$\ad_{L_E}(\overline  F_1) = \r \overline F_1\ ,\qquad \ad_{L_E}(\overline F_2) = \r \overline F_2 + \overline F_1$$
where $\r$ is one of the   two possibilities  given in \eqref{possible-eigen}.
Plugging this into \eqref{res} with $\overline F = \overline  F_2$,  we obtain that
$$ \r^2 \overline F_2 + 2 \r \overline F_1 = - \a(E) \r \frac{ 2 \mu -3}{2} \overline F_2  - \a(E)  \frac{2 \mu -3}{2} \overline F_1 - \a(E)^2 \frac{(\mu-1)(\mu -2)}{4} \overline F_2\ .
$$
Looking at the coefficients of $\overline F_1$, we see that $\r$ and $\a(E)$ have to satisfy the relation
$$4 \r =  - \a(E) (2\mu -3)\ .$$
Since $\a(E) \neq 0$ and $4 \r$ is either $-\a(E)(2\mu-2)$ or $-\a(E)(2\mu-4)$, we get a contradiction in both cases.  This  shows that only (a) and (b) can occur.
\par
\medskip
To conclude the proof, it remains only  to show  that $\ad_{L_E}|_{\gm^{01}}$  is diagonalisable  also in case $\a(E) = 0$.
Suppose not. Then, considering Jordan canonical forms as before,  we may pick a basis $(\overline F_1, \ldots, \overline F_n)$ for $\gm^{01}$  such that
\setcounter{section}{5}
\setcounter{equation}{17}
\beq  \label{pippo-2} L_{E|\overline F_1}  =0\ ,\qquad L_{E|\overline F_2}  = \overline F_1\ .\eeq
We may also consider a basis $(E_1, \ldots, E_n)$ for $\gm^{10}$  with  $E_1 = E$ and the corresponding  basis $\cB = (E_i, \overline F_j)$ for $V^\bC = (\gm^{-1})^\bC$.
We  start with  a couple of  preliminary observations. We first notice  that,  by the
 Jacobi identity and  by the property   $[\gm^{01}, \gm^{01}] = \{0\}$, for any $1 \leq i,j \leq n$
\beq \label{pippo-1}\overline F_i \bdot E \bdot \overline F_j = [\overline F_i, [E, \overline F_j]] = [[\overline F_i, E], \overline F_j] = \overline F_j \bdot E \bdot \overline F_i\ .\eeq
Second, by  $L_{E|E} =\a(E) E = 0$ and  \eqref{war1-CRbis},  we have that   for   $1 \leq i, j \leq 2$ and   $1 \leq r$
\beq
\begin{split}\label{pippo} & L_{E^r \bdot (\overline F_i \bdot E \bdot \overline F_j)}  =
r E^{r-1}\bdot   L_{E| \overline F_i \bdot E \bdot \overline F_j} -  E^{r-1}\bdot   L_{\overline F_i \bdot E \bdot \overline F_j |E}  = \\
&= r E^{r-1}\bdot   L_{E| \overline F_i}  \bdot E \bdot \overline F_j + r E^{r-1}\bdot   \overline F_i  \bdot E \bdot L_{E| \overline F_j}
-   E^{r-1}\bdot   L_{\overline F_i \bdot E \bdot \overline F_j |E}  \overset{\eqref{pippo-2}}= \\
&= r E^{r-1}\bdot \left( \d_{i2} \overline F_1 \right)\bdot E \bdot \overline F_j + r E^{r-1}\bdot  \overline F_i  \bdot E \bdot \left(\d_{j2} \overline F_1 \right)  -   E^{r-1}\bdot   L_{\overline F_i \bdot E \bdot \overline F_j |E} \overset{\eqref{pippo-1}}= \\
&=  r E^{r-1}\bdot  (\d_{i2} \overline F_j    + \d_{j 2} \overline F_i )\bdot E \bdot \overline F_1 -    E^{r-1}\bdot   L_{\overline F_i \bdot E \bdot \overline F_j |E}\ .
\end{split}
\eeq
We now observe that the term $ L_{\overline F_i \bdot  E \bdot \overline F_j }$ is in $[\gg^1, (\gm^{-3})^\bC] \subset (\gm^{-2})^\bC$. Thus it has the  form
$L_{\overline F_i \bdot  E \bdot \overline F_j }   = \sum_{k, \ell} c^{\ell k}_{(ij)} E_\ell\bdot \overline F_k$
 for some appropriate   coefficients $c^{\ell k}_{(ij)}$ and
  \beq
\label{pippo-bis}  L_{E^r \bdot (\overline F_i \bdot E \bdot \overline F_j)}  =  r E^{r-1}\bdot  (\d_{i2} \overline F_j    + \d_{j 2} \overline F_i )\bdot E \bdot \overline F_1 -    \sum_{k, \ell} c^{\ell k}_{(ij)} E^{r-1}\bdot   E_\ell\bdot \overline F_k \bdot E\ .
\eeq
 We are now ready to conclude.  Since $E^{\mu-2} \bdot \overline F_2 \bdot E \bdot \overline F_2 = 0$ (it has degree $\mu + 1$),
\begin{multline} 0 = L_{E^{\mu-2} \bdot (\overline F_2 \bdot E \bdot \overline F_2)} \overset{\eqref{pippo-bis}} = 2 (\mu-2) E^{\mu -3}\bdot \overline F_2 \bdot E \bdot \overline F_1 -     \sum_{k, \ell} c^{\ell k}_{(22)} E^{\mu-3}\bdot   E_\ell\bdot \overline F_k \bdot E
\ ,\end{multline}
which implies the  following relation  between elements of degree $\mu$
$$E^{\mu -3}\bdot \overline F_2 \bdot E \bdot \overline F_1  = \frac{1}{2(\mu -2)} \sum_{k, \ell} c^{\ell k}_{(22)} E^{\mu-3}\bdot   E_\ell\bdot \overline F_k \bdot E\ .$$
If we now apply $\ad_L$ to both sides  and use \eqref{pippo-bis} once again, we obtain the following identity  between homogeneous elements of degree $\mu -1$:
\begin{multline} \label{pippo-ter}  E^{\mu -4}\bdot \overline F_1 \bdot E \bdot \overline F_1 = \\
 =  \frac{1}{\mu -3} \sum_{k, \ell} c^{\ell k}_{(21)} E^{\mu -4}\bdot   E_\ell\bdot \overline F_k \bdot E
 +
\frac{1}{2(\mu -2)(\mu -3)}  \sum_{k, \ell} c^{\ell k}_{(22)} L_{E^{\mu-3}\bdot   E_\ell\bdot \overline F_k \bdot E}\ .
\end{multline}
We now observe  that the left hand side  is an element in $\gm^\bC = \gU_{V^\bC}/(\gi_{10} + \gi_{01} + \gi^{-\mu})$, which is the projection of a  non-zero $\cB$-monomial in $\gU_{V^\bC}$ whose $\cB$-type is
\beq \label{ecco-bis} (\underset{\text{entries  corresp. to  the } E_i}{\underbrace{\mu- 3,0, \ \ldots\ , 0}},\underset{\text{entries  corresp. to  the } \overline F_k}{\underbrace{2, 0,0, \ \ldots\ , 0}})\ .\eeq
On the other hand, a straightforward check shows  that the right hand side  is a linear combination of  projections of $\cB$-monomials,   whose  $\cB$-types might only have   the form
$$(\underset{\text{entries  corresp. to  the } E_i}{\underbrace{\ast ,\ \ \ast,\ \  \ldots, \ \ \ast,\ \ \ast}}, \underset{\text{entries  corresp. to  the } \overline F_k}{\underbrace{0, \ldots, 0,  1, 0, \ldots, 0}})\ .$$
 Since these $\cB$-types  are surely different from \eqref{ecco-bis},
  Corollary \ref{cor44} implies that the left hand side  and the right hand side of   \eqref{pippo-ter}  are linearly independent and that  the  equality cannot be true. This contradiction concludes the proof. \end{pf}
\begin{cor} \label{corollary} Given  $L \in \gg^1$,  with associated linear map $\a = \a_L \in (\gm^{10})^*$, there exists a basis $(E_1, \ldots, E_n, \overline F_1, \ldots, \overline F_n)$ for $\gm^\bC$, with $E_i \in \gm^{10}$, $\overline F_j \in \gm^{01}$,  and an  $n$-tuple of rational numbers $\r_i$, $1 \leq  i  \leq  n$, each of them equal either  to $-\frac{\mu - 1}{2}$ or $-\frac{\mu-2}{2}$,    such that $E_1, \ldots, E_{n-1} \in \ker \a$ and  for any $E \in \gm^{10}$ and $ 1 \leq i, k \leq n$
\beq \label{5.17}  L_{E|\overline F_i} = -  \a(E) \r_i \overline F_i\ , \qquad L_{E|E_k}  =   \frac{\d_{kn} \a(E_n)}{2}  E +  \frac{ \a(E)}{2} E_k\ .\eeq
In particular, $L  \neq  0$  if and only if $\a \neq 0$.
\end{cor}
\begin{pf}  By Proposition \ref{cor1}, the claim is trivial if $\a = 0$.  Then, we may assume that   $\a \in \Hom(\gm^{10}, \bC)$ is non-trivial and we consider  a basis $(E_1, \ldots, E_n)$ for $\gm^{10}$ such that $(E_1, \ldots, E_{n-1})$ is a basis for $\ker \a$ and  $E_n \notin \ker \a$. By Proposition \ref{cor1}, the restricted adjoint action $\ad_{L_{E_n}}|_{\gm^{01}}$ is diagonalisable and there exists a basis $(\overline F_1, \ldots, \overline F_n)$ for $\gm^{01}$, made of eigenvectors whose  associated eigenvalues have the form  $- \a(E_n)\frac{\mu - 1}{2}$ or $-\a(E_n)\frac{\mu-2}{2}$. The first identity in  \eqref{5.17} follows by the facts that, for an arbitrary  $E = \l^k E_k  \in \gm^{10}$, one has that
$$\a(E) = \sum_{k= 1}^{n-1} \l^k \a(E_k) + \l^n \a(E_n) = \l^n \a(E_n)$$
 and that, setting $\rho_i$ as in the statement, for each $\overline F_i$
$$L_{E|\overline F_i} =  \sum_{k = 1}^{n-1} \l^k L_{E_k|\overline F_i} +  \l^n L_{E_n|\overline F_i} = - \l^n  \a(E_n) \r_i \overline F_i =  -  \a(E) \r_i \overline F_i\ .$$
The second identity of \eqref{5.17} comes from \eqref{5.13}. The last claim is a   consequence of \eqref{5.17} and  the fact  that $L \neq 0$ if and only if there exists some  $E \in \gm^{10}$  for which the element  $L_E \in (\gg^0)^\bC$ is non trivial.
\end{pf}
\subsection{Proof of Theorem \ref{main}}
Let $\gm + \gg^0$ be a  totally nondegenerate CR Tanaka  algebra of depth $\mu\geq 4$ and $\gm^\bC + (\gg^0)^\bC = \gU_J/\gi^{-\mu} + (\gg^0)^\bC$ be its complexification. To prove the theorem, it suffices  to show  that a  contradiction arises if one assumes that there is  $L  \neq 0$ in  $\gg^1$ (which, as usual,  we   identify with the corresponding element in $\jmath(\gg^1) \subset \wh \gg^1$ of the prolongation of $\gm^\bC + (\gg^0)^\bC$).
\par
\smallskip
So, let us assume  the existence of  $0 \neq L \in \gg^1$ with associated non-trivial linear map $\a = \a_L \neq 0$. We may therefore consider a basis $\cB = (E_k, \overline F_i)$ with the properties described in Corollary \ref{corollary} and with $\a(E_n) \neq 0$.   We may also assume that $E_n = F_k$ for some $k$. In fact,  the considered   basis $(\overline F_1, \ldots, \overline F_n)$ for $\gm^{01}$ can be actually chosen
 in such a way  to be  independent of $E$, a fact that can  be checked as follows.  Pick an element $E_o \in \gm^{10}$ such that $\a(E_o) = 1$ and  for any  other  $E \in \gm^{10}$,  set $c_E \= \a(E)$. Since
$ E - c_E E_o \in \ker \a$, by Proposition \ref{cor1}, it follows that    $\ad_{L_{E- c_E E_o}}|_{\gm^{01}} = 0$.
This  implies that $\ad_{L_{E}}|_{\gm^{01}} = c_E \ad_{L_{E_o}}|_{\gm^{01}}$ and  that  a diagonalising basis $(\overline F_1, \ldots, \overline F_n)$ for $\ad_{L_{E_o}}|_{\gm^{01}}$ is necessarily a  diagonalising basis for   $\ad_{L_{E}}|_{\gm^{01}}$ for any possible choice of $E$. Now, being such a basis fixed,  for at least one  element $\overline F_k$ we must have   $\a(F_k) \neq 0$  (otherwise,  we would have that $L|_{\gm^{10}}   = 0$ and this  would imply $L = 0$). We may therefore
 assume that  $E_n  = F_k$ as claimed. \par
Hence, by possibly reordering the basis  $(\overline F_i)$ for $\gm^{01}$ and by an appropriate rescaling, there is no loss of generality if we have  $E_n = F_n$ and $\a(E_n) = 1$. \par
\smallskip
In what follows,  for simplicity of notation, we set $E \= E_n$, $\overline E \= \overline F_n$ and $ \r \= \r_n$.  \par
\smallskip
Our proof is crucially based on  the following identity, which holds  for each $r \geq 0$:
\beq \label{claim} L_ {E^r \bdot  (\overline E \bdot (E\bdot \overline E))} =  \frac{r(4 \r + r + 1)}{2}  E^{r-1} \bdot (\overline E \bdot (E \bdot \overline E)) + ( 2 \rho  + 1) E^{r+1}\bdot \overline E \ .
\eeq
In fact, recalling that $[E, E] = 0 = [\overline E, \overline E]$  and that  $L_{E|E} = E$, $L_{E|\overline E} = \r \overline E$,  $L_{\overline E|E} = \overline{L_{E|\overline E}} = \r  E$, the identity \eqref{claim}  for  $r = 0$ holds because
\begin{multline}\label{5.19} L_{\overline E \bdot (E \bdot \overline E)} = [L, [\overline E, [E, \overline E]]] =  [[L, \overline E], [E, \overline E]]  + [\overline E, [L, [E, \overline E]]]= \\
 =  [[[L, \overline E], E], \overline E]  +  [E,[ [L, \overline E],  \overline E] ] +  [\overline E,[ [L, E], \overline E]]] +  [\overline E,  [E, [L,\overline E]]] = \\
 = L_{\overline E|E} \bdot \overline E + E \bdot L_{\overline E| \overline E} -L_{E|\overline E}  \bdot  \overline E + L_{\overline E| E} \bdot  \overline E  = (2 \rho +1) E \bdot \overline E \ .
\end{multline}
When $r = 1$,  the identity follows from  the fact that $L_E \in \gg^0$  (which  means that $L_E$ is a derivation of the Lie algebra $\gm^\bC$) and from  \eqref{5.19}. Indeed, all this  implies
 \begin{multline}\label{5.20} L_{E\bdot (\overline E \bdot (E \bdot \overline E))}  =[ L_E, \overline E \bdot (E \bdot \overline E)] + [E, L_{\overline E \bdot (E \bdot \overline E)}] = \\
  = L_{E|\overline E} \bdot (E \bdot \overline E) + \overline E\bdot (L_{E|E} \bdot \overline E) + \overline E \bdot (E \bdot L_{E|\overline E}) + (2 \rho +1) E^2\bdot \overline E =\\
  =( 2 \rho  + 1) (\overline E \bdot (E \bdot \overline E) + E^2\bdot \overline E) \ .
 \end{multline}
 Finally, when $r \geq 2$,   by \eqref{war1-CRbis},
  \begin{multline}
  L_ {E^{r} \bdot  (\overline E \bdot (E\bdot \overline E))}  = \\
  = r E^{r-1} \bdot L_{E|\overline E \bdot (E\bdot \overline E)}  -   E^{r-1} \bdot   L_{\overline E \bdot (E\bdot \overline E)|E}  +  \frac{r (r-1)}{2}
E^{r-2}\bdot (L_{E|E} \bdot   (\overline E \bdot (E\bdot \overline E))) = \\
= r E^{r-1} \bdot (L_{E|\overline E} \bdot (E\bdot \overline E))  + r E^{r-1} \bdot (\overline E \bdot (L_{E|E}\bdot \overline E)) + r E^{r-1} \bdot (\overline E \bdot (E\bdot L_{E| \overline E}))  +\\
 +  E^{r} \bdot   L_{\overline E \bdot (E\bdot \overline E)}  + \frac{r (r-1)}{2}
E^{r-1}\bdot  (\overline E \bdot (E\bdot \overline E)) = \\
 = (2 \r +1) r E^{r-1} \bdot (\overline E \bdot (E\bdot \overline E)) +  (2 \r +1) E^{r +1} \bdot  \overline E +  \frac{r (r-1)}{2}
E^{r-1}\bdot  (\overline E \bdot (E\bdot \overline E)) \ ,
\end{multline}
  which gives \eqref{claim}  in all remaining cases.\par
\medskip
Let  us now consider  the first integer  $\nu$ such that $E^{\nu}\cdot (\overline E \cdot (E\cdot \overline E))= 0$.  Since  $\gm + \gg^0$ is  totally nondegenerate, then  $\nu = \mu-3$ or $\mu -2$.  We  consider these two cases separately and show that in both cases   the desired contradiction arises.\par
 \medskip

 \noindent
 {\it Case $1$: $\nu=\mu-3$.} \\
  Since  $E^{\mu-3}\bdot(\overline E\bdot(E\bdot\overline E))) = 0$, from  \eqref{claim} we have
\begin{multline*}
0=L_{E^{\mu-3}\bdot(\overline E\bdot(E\bdot\overline E)))}=\\
= \frac{(\mu-3)(4\rho+\mu-2)}{2}E^{\mu-4}\bdot(\overline E\bdot(E\bdot\overline E))+(2\rho+1)E^{\mu-2}\bdot\overline E\ .
\end{multline*}
The  two   monomials in the right hand side are in $(\gm^{-(\mu-1)})^\bC=\gU^{-(\mu-1)}_J$ and have different  $\cB$-types. Hence, they are linearly independent and since $\mu \geq 4$, this contradicts the fact that, for both possibilities  $\rho= - \frac{\mu-1}{2}$ and $\r = - \frac{\mu-2}{2}$,  the coefficients of the  linear combination in the right hand side   are  non-zero.
\par
 \medskip
 \noindent
 {\it Case $2$: $\nu=\mu-2$}. \\
  From  $E^{\mu-2}\bdot(\overline E\bdot( E\bdot\overline E)) = 0$ and \eqref{claim}
  \begin{multline*}
 0=L_{E^{\mu-2}\bdot(\overline E\bdot(E\bdot\overline E))}=\\
 = \frac{(\mu-2)(4\rho+\mu-1)}{2}E^{\mu-3}\bdot(\overline E\bdot(E\bdot\overline E))+(2\rho+1)E^{\mu-1}\bdot\overline E\ .
 \end{multline*}
 This time the right hand side is an element in $(\gm^{-\mu})^\bC=\gU^{-\mu}_J/\gi^{-\mu}$
 and we cannot claim that the two monomials appearing in such expression are linearly independent.  However, if we apply $L$ to this expression and use \eqref{claim} and \eqref{war1-CRbis} we get
\begin{multline*}
  0=\frac{(\mu-2)(4\rho+\mu-1)}{2}\left(\frac{(\mu-3)(4\rho+\mu-2)}{2}E^{\mu-4}\bdot(\overline E\bdot( E\bdot\overline E))+\right.\\
\hskip 8cm   + (2\rho+1)E^{\mu-2}\bdot\overline E\bigg) +
  \\
  +(2\rho+1)\left(\rho(\mu-1)E^{\mu-2}\bdot\overline E+\frac{(\mu-1)(\mu-2)}{2}E^{\mu-2}\bdot\overline E\right) = \\
  =  \frac{(\mu-2)(\mu-3) (4\rho+\mu-1)(4\rho+\mu-2)}{2}E^{\mu-4}\bdot(\overline E\bdot( E\bdot\overline E)) + \\
  + c E^{\mu-2}\bdot\overline E
\end{multline*}
with $c =  (2\rho+1) \left(\frac{(\mu-2)(4\rho+\mu-1) }{2} + \rho(\mu-1) +  \frac{(\mu-1)(\mu-2)}{2}\right)$.
Now,  both   monomials in the right hand side are in $(\gm^{-(\mu-1)})^\bC=\gU^{-(\mu-1)}_J$ and have different  $\cB$-types. Thus they are linearly independent and  the coefficient $\frac{(\mu-2)(\mu-3) (4\rho+\mu-1)(4\rho+\mu-2)}{2}$ of $E^{\mu-4}\bdot(\overline E\bdot( E\bdot\overline E))$ should be  $0$. This is in contradiction with the fact that, being $\mu \geq 4$, such coefficient is nonzero for both possibilities     $ \r = - \frac{\mu-1}{2}$ and $\r = - \frac{\mu-2}{2}$.
 \par

\end{document}